\documentclass{article}

\usepackage[utf8]{inputenc}
\usepackage{amsmath}
\usepackage{fullpage}
\usepackage{relsize}
\usepackage{amssymb}
\usepackage{amsfonts}
\usepackage{graphicx}
\usepackage{tikz}
\usepackage{float}
\usepackage{caption}
\usepackage{authblk}
\usepackage[english]{babel}
\usepackage{mathtools}
\usepackage{pifont}
\usepackage{subfig}
\usepackage{graphicx}
\usepackage{float}
\usepackage{xcolor,colortbl}

\makeatletter
\newcommand{\thickhline}{%
    \noalign {\ifnum 0=`}\fi \hrule height 1.5pt
    \futurelet \reserved@a \@xhline
}
\newcolumntype{"}{@{\hskip\tabcolsep\vrule width 1.5pt\hskip\tabcolsep}}
\newcolumntype{?}{!{\vrule width 1.5pt}}
\makeatother

\usepackage[square,comma,numbers,sort&compress]{natbib}
\usepackage{hyperref}

\usetikzlibrary{arrows.meta}
\usetikzlibrary{calc}
\usetikzlibrary{decorations.pathmorphing, patterns,shapes}

\renewcommand{\i}{\mathrm{i}}
\newcommand{\e}{\mathrm{e}}
\newcommand{\eps}{\epsilon}

\newcommand{\norm}[1]{|| #1 ||}

\renewcommand{\d}{\,\mathrm{d}}
\newcommand{\diff}[2]{\frac{\mathrm{d} #1}{\mathrm{d} #2}}

\usepackage{todonotes}

\title{Exponential Asymptotics using Numerical Rational Approximation in Linear Differential Equations}
\author[1,3]{Christopher J. Lustri\footnote{Electronic address: christopher.lustri@sydney.edu.au}}
\author[2,3]{Samuel Crew}
\author[4]{S. Jonathan Chapman}
\affil[1]{School of Mathematics and Statistics, The University of Sydney, New South Wales, 2006, Australia}
\affil[2]{Faculty of Computer Science, Ruhr Universit\"{a}t Bochum,
Universit\"{a}tsstra{\ss}e, Bochum, 44799, Germany}
\affil[3]{Theoretical Sciences Visiting Program, Okinawa Institute of Science and Technology Graduate University, Onna, 904-0495, Japan.}
\affil[4]{Mathematical Institute, University of Oxford, Oxford, OX2 6GG, United Kingdom}
\date{}

\begin{document}

\maketitle

\abstract{Singularly-perturbed ordinary differential equations often exhibit Stokes' phenomenon, which describes the appearance and disappearance of oscillating exponentially small terms across curves in the complex plane known as Stokes curves. These curves originate at singular points in the leading-order solution to the differential equation. In many important problems, it is impossible to obtain a closed-form expression for these leading-order solutions, and it is therefore challenging to locate these singular points. We present evidence that the analytic leading-order solution of a linear differential equation can be replaced with a  rational approximation based on a numerical leading-order solution using the adaptive Antoulas-Anderson (AAA) method. We show that the subsequent exponential asymptotic analysis accurately predicts the exponentially small behaviour present in the solution. We explore the limitations of this approach, and show that for sufficiently small values of the asymptotic parameter, this approach breaks down; however, the range of validity may be extended by increasing the number of poles in the rational approximation. We finish by presenting a related nonlinear problem and discussing the challenges that arise when attempting to apply this method to nonlinear problems.}

\section{Introduction}

The behaviour known as ``Stokes' phenomenon'' was first observed by George Gabriel Stokes in his analysis of solutions to the Airy equation \cite{Stokes,stokes1902discontinuity}. The Airy function displays oscillatory behaviour for large negative arguments, and decays exponentially for large positive arguments. Stokes wished to understand the transition between these two regimes by studying the large-argument asymptotic behaviour in the complex plane. In doing so, Stokes discovered that exponentially-small correction terms in the asymptotic expansion appear or disappear as special curves in the complex plane are crossed. These curves are known as ‘Stokes curves’.

New asymptotic techniques were developed in order to study the switch in behaviour across Stokes curves. Berry \cite{Berry1988,Berry,Berry1991}, and Berry \& Howls \cite{BerryHowls1990,berry1991hyperasymptotics} showed that truncating a a divergent asymptotic series optimally results in asymptotic error of exponentially small order. Hence, rescaling the problem to study this truncation error directly allows for the calculation of these exponentially small terms. This idea was used to develop new asymptotic methods for studying exponentially small behaviour, known as exponential asymptotics, asymptotics-beyond-all-orders or hyperasymptotics \cite{ berry1991hyperasymptotics, daalhuis1993hyperasymptotics}. We will denote exponentially small contributions to the asymptotic approximation of some solution $u$ as $u_{\mathrm{exp}}$.

Exponential asymptotic techniques typically require the explicit calculation of the leading-order solution of the differential equation, which we will denote as $u_0$.  Stokes curves originate at singularities of $u_0$ \cite{Dingle}. Additionally, the asymptotic size of $u_{\mathrm{exp}}$ as $\eps \to 0$ typically depends on the location and strength of the singularity. For many problems, knowing the analytical form of the leading-order solution near singularities provides enough information to calculate the Stokes curves and approximate the exponential contributions \cite{Chapman}. 

For many practical problems, it is impossible to calculate the leading-order solution analytically \cite{Deng2021,Deng2022,deng2023exponential,ChapmanTrinh}. In this case, the leading-order solution must instead be approximated numerically. There is growing evidence that numerical methods for analytic continuation, including numerically stepping into the complex plane \cite{ChapmanTrinh} and rational approximation methods \cite{deng2023exponential} can be used in conjunction with exponential asymptotics. 

In general, the leading-order solution to an ordinary differential equation can be computed on a real domain using any one of a variety of standard numerical methods \cite{hildebrand1987introduction}. Using the result of numerical computations to describe the complex plane behaviour of these solutions is challenging. It is well-known that analytic continuation is ill-posed, and that small changes to the function being continued, such as numerical error associated with computational methods, can lead to significant changes in the analytically-continued result. 

Despite this, significant progress has been made on performing analytic continuation using numerical methods; this is typically performed by requiring the output function to take a particular algebraic form. The most commonly-used methods for numerical complex analytic continuation are methods that give a rational function as an output. These rational approximation methods include Pad\'{e} approximation \cite{baker1961pade} and the adaptive Antoulas--Anderson (AAA) algorithm \cite{Nakatsukasa}. 

 We wish to consider settings in which the solution has been numerically approximated on some discrete set of points along the real axis which we then seek to analytically continue into the complex plane. This setup makes the AAA algorithm particularly well-suited for our purpose. This algorithm is described in Section \ref{S.NAC}, and produces a rational function, denoted here as $\hat{u}_0$ which approximates the leading-order behaviour. Note that all singularities in $\hat{u}_0$ must be simple poles; if $u_0$ contains other singularities, such as branch points, the branch will be typically be approximated as an accumulation of simple poles instead \cite{Gopal,Trefethen2020,Trefethen2021}. Recall that the strength of the singularity in $u_0$ is generically different to the simple poles present in $\hat{u}_0$, and the asymptotic behaviour of $u_{\mathrm{exp}}$ depends on the strength of the singularity. This argument indicates that any exponential term obtained using this leading-order approximation must display different asymptotic behaviour to ${u}_{\mathrm{exp}}$ in the limit that $\eps \to 0$.
 
In this study, we focus our attention on a singularly-perturbed linear ordinary differential equation in the limit that $\eps \to 0$, whose leading-order solution $u_0$ contains two branch points, namely 
\begin{equation}\label{eq.ode}
\epsilon^2 u''(x) + u(x) = \frac{1}{\sqrt{x+\i}} + \frac{1}{\sqrt{x-\i}}, 
\end{equation}
with the boundary conditions
\begin{equation}\label{eq.bc}
\lim_{x\to-\infty} u(x) = 0, \qquad \lim_{x\to-\infty} u'(x) = 0.
\end{equation}
We note that differential equation \eqref{eq.ode} possesses a closed-form solution in terms of the exponential integral function,
\begin{align}\nonumber
    u(x) =  \,&A(\eps) \cos\left(\frac{x}{\epsilon}\right) + B(\eps) \sin\left(\frac{x}{\epsilon}\right) \\
    &- \frac{\i\sqrt{x+\i}}{2\eps}\e^{(1-\i x)/\eps}\int_1^{\infty} \frac{\e^{-(1-\i x)t/\eps}}{\sqrt{t}} \d t + \frac{\i\sqrt{x+\i}}{2\eps}\e^{-(1-\i x)/\eps}\int_1^{\infty} \frac{\e^{(1-\i x)t/\eps}}{\sqrt{t}} \d t + \mathrm{c.c.},
    \end{align}
where $A$ and $B$ may be determined from the boundary conditions. It is possible to determine $A$ and $B$ by performing a careful analysis on the integral expressions in the limit that $x \to -\infty$, but we will instead take a different approach, and determine an asymptotic expansion for the solution in the limit that $\eps \to 0$.

We will perform a typical exponential asymptotic analysis on \eqref{eq.ode}--\eqref{eq.bc}, using the method proposed in \cite{Daalhuis}, to determine ${u}_{\mathrm{exp}}$ for this problem. We will then sample $u_0$ on a discrete set of points and use this as the input for a rational approximation $\hat{u}_0$. We will apply the same exponential asymptotic method on this leading-order expression to obtain $\hat{u}_{\mathrm{exp}}$. By comparing the analytic form of ${u}_{\mathrm{exp}}$ and $\hat{u}_{\mathrm{exp}}$, we will see that the two expressions have different asymptotic behaviour in the limit $\epsilon \to 0$. 

We will then show that, despite this difference, $\hat{u}_{\mathrm{exp}}$ is able to accurately approximate $u_{\mathrm{exp}}$ for a range of values of $\epsilon$, and that this range can be extended by increasing the accuracy of the AAA approximation (ie. reducing the $L^2$ error threshold on the set of sample points). Finally, we will provide evidence that the threshold value of $\epsilon$ below which $\hat{u}_{\mathrm{exp}}$ is inaccurate is proportional to the difference between the true branch point location in $u_0$ and the pole nearest to this point in $\hat{u}_0$. This may be thought of as the approximation error of the true branch point location from the rational approximation algorithm.

\subsection{Numerical Analytic Continuation}\label{S.NAC}

Rational approximation methods typically approximate some function $f(x)$ as a rational expression,
\begin{align}
f(x) \approx\frac{n(x)}{d(x)},
\end{align}
where $n(x)$ and $d(x)$ are polynomials.

Historically, the most widely-used method for obtaining rational approximations is Pad\'{e} approximation, which is constructed using the Taylor coefficients of $f(x)$ at a particular point. We want to study data that takes a different form; where values of $f(x)$ are calculated numerically on some discrete support set of points. While it is possible to approximate Taylor coefficients using values of the function on a discrete set of points, there are other methods more well-suited to input data with this form.

On such is the AAA approximation, which works by iteratively solving a minimisation problem to produce an optimal rational approximation. The explanation that follows is a high-level description of the AAA algorithm; for a detailed explanation, see \cite{Nakatsukasa}. 

The algorithm is based on expressing the rational approximation in the form
\begin{align}
f(x) \approx\frac{n(x)}{d(x)}={\sum_{j=1}^m\frac{w_jf_j}{x-x_j}} \, \Bigg/ \, {\sum_{j=1}^m\frac{w_j}{x-x_j}},
\label{e:intro_aaa}
\end{align}
The points $x_j$ for $j = 1, 2, \ldots, m$, are known as support points. They are drawn from the sample set $X$, which consists of the points on which $f(x)$ is known; we define $f_j=f(x_j)$.  The algorithm solves an optimization problem in order to determine the weights $w_j$, and then generates another support point $x_{j+1}$. This process then repeats iteratively.

The weights $w_j$ are generated by solving a linear least-squares problem over a set of points in the restricted domain $X^{(m)}=X \setminus \{x_1,\ldots,x_{m}\}$; these points are labelled $X_i^{(m)}$. The weight vector $w =(w_1,w_2,\dots,w_m)^T$ is chosen to minimize $\norm{fd-n}_{X^{(m)}}$ subject to $\norm{w}_m=1$, where $\norm{\cdot}_{X^{(m)}}$ is the discrete 2-norm over $X^{(m)}$ and $\norm{\cdot}_m$. is the discrete $L^{2}$ on $m$-vectors. If the 2-norm of the approximation residuals on $X^{(m)}$ are beneath some specified tolerance (relative to the maximum value taken by $|f(x_j)|$), the algorithm terminates.

If the algorithm does not terminate, the value of $x_{m+1}$ is determined for the next iteration. It is found by choosing the value of $x_{m+1}\in X^{(m)}$ that maximizes the quantity
\begin{align}
\sum_{j=1}^m \frac{w_j f(x_{m+1})}{x_{m+1}-x_j}-
\sum_{j=1}^m \frac{w_j f_j}{x_{m+1}-x_j}.
\end{align}
The process then repeats until the termination criteria is reached.

The key information for our purposes is that the method takes in the function values on a set of points, and returns a rational approximation that minimizes the $L^2$-norm of the approximation error on the data support set, and that it is an iterative process that increases the number of poles in the solution until some predetermined tolerance is met. We often write the approximation in the form
\begin{equation}\label{e.AAA}
f(x) \approx \sum_{j=0}^{m} \frac{a_j}{x - p_j},
\end{equation}
where $a_j$ are the residues of the function at each pole location $p_j$. 

As previously alluded to, one potential obstacle in applying rational approximation methods is that these methods produce meromorphic functions; all singularities in the approximation are isolated simple poles. Previous studies of AAA approximation~\cite{Gopal,Trefethen2020,Trefethen2021} showed that the rational function approximation of a target function with a branch point contains an exponential clustering of poles (and zeroes) in the approximation approaching the branch point. This configuration of poles accurately approximates the effects of the branch cut. The distribution of poles when approximating branch cuts in rational approximations has been studied rigourously in the related problem of Pad\'{e} approximation \cite{stahl1997convergence}.

In the following we show that if we approximate the leading-order solution of an ordinary differential equation with a rational function of the form \eqref{e.AAA}, it is possible to use exponential asymptotics on this rational function in order to approximate the exponentially small terms that appear in the solution due to Stokes' phenomenon.

\section{Exponential Asymptotic Analysis}\label{S.expasymp}

\subsection{Asymptotic Power Series}

\begin{figure}[tb]
\centering
 \begin{tikzpicture}
 [x=15,>=stealth,y=15]

\draw[blue, only marks,mark=*] plot[] file {Fig_Div_Data.txt};
\draw (-0.5,-5) node[left] {\scriptsize{$-5$}} -- (-0.5,3) node[left] {\scriptsize{$3$}} -- (16.5,3) -- (16.5,-5) -- cycle;
\draw[dotted] (-0.5,-4) node[left] {\scriptsize{$-4$}} -- (16.5,-4);
\draw[dotted] (-0.5,-3) node[left] {\scriptsize{$-3$}} -- (16.5,-3);
\draw[dotted] (-0.5,-2) node[left] {\scriptsize{$-2$}} -- (16.5,-2);
\draw[dotted] (-0.5,-1) node[left] {\scriptsize{$-1$}} -- (16.5,-1);
\draw[->] (-0.5,0) node[left] {\scriptsize{$0$}} -- (17,0) node[right] {\scriptsize{$n$}};
\draw[dotted] (-0.5,1) node[left] {\scriptsize{$1$}} -- (16.5,1);
\draw[dotted] (-0.5,2) node[left] {\scriptsize{$2$}} -- (16.5,2);

\node[rotate=0,above] at (7.5,3) {\scriptsize{$\log_{10}|\eps^{2n} u_n|$ at $x = 0$}};

\draw (0,-5) node[below] {\scriptsize{$0$}} -- (0,-4.7);
\draw (1,-5) -- (1,-4.85);
\draw (2,-5) -- (2,-4.85);
\draw (3,-5) -- (3,-4.85);
\draw (4,-5) -- (4,-4.85);
\draw (5,-5) node[below] {\scriptsize{$5$}}-- (5,-4.7);
\draw (6,-5) -- (6,-4.85);
\draw (7,-5) -- (7,-4.85);
\draw (8,-5) -- (8,-4.85);;
\draw (9,-5) -- (9,-4.85);
\draw (10,-5) node[below] {\scriptsize{$10$}}-- (10,-4.7);
\draw (11,-5) -- (11,-4.85);
\draw (12,-5) -- (12,-4.85);
\draw (13,-5) -- (13,-4.85);
\draw (14,-5) -- (14,-4.85);
\draw (15,-5) node[below] {\scriptsize{$15$}}-- (15,-4.7);
\draw (16,-5) -- (16,-4.85);

\draw[thick,blue,<-] (5,-3.9) -- (5,-2.8) node[above] {\scriptsize{Minimum value}};

 \end{tikzpicture}
\caption{
Magnitude of series terms from \eqref{eq.series} evaluated at $x = 0$ for $\epsilon = 0.1$. As $n$ increases, the terms become smaller until a minimum value is reached at $n = 5$, after which the terms increase in size due to the factorial contribution to the numerator of $u_n$ in \eqref{eq.seriesterms} The series \eqref{eq.series} must therefore be divergent, and the optimal truncation point occurs at the minimum value.}\label{Fig_seriesterms}
\end{figure}
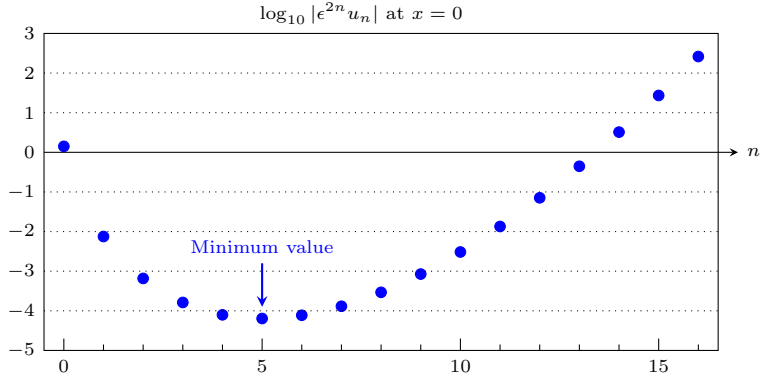

We propose an asymptotic series solution to \eqref{eq.ode}--\eqref{eq.bc} of the form
\begin{equation}\label{eq.series}
u(x; \epsilon) \sim \sum_{n=0}^{\infty} \epsilon^{2n}u_n(x).
\end{equation}
By substituting this into the ordinary differential equation \eqref{eq.ode} and matching powers of $\epsilon$ in the limit that $\epsilon \to 0$, we find that the leading-order behaviour is given by
\begin{align}\label{eq.LO}
u_0 = \frac{1}{\sqrt{x+\i}} + \frac{1}{\sqrt{x-\i}}.
\end{align}
Matching at $\mathcal{O}(\eps^n)$ as $\eps \to 0$ gives
\begin{align}\label{eq.recur}
u_n = -u_{n-1}'', \qquad n \geq 1.
\end{align}
This recurrence relation can be used calculate $u_n$ exactly, giving
\begin{equation}\label{eq.seriesterms}
u_n = \frac{(4n)!(-1)^n}{2^{4n}(2n)!(x-\i)^{2n+{1/2}}} + \mathrm{c. c.},
\end{equation}
The series \eqref{eq.series} with terms given by \eqref{eq.seriesterms} is divergent for any fixed choice of $x$. Figure \ref{Fig_seriesterms} depicts $\eps^n u_n$ at $x = 0$ with $\eps = 0.1$. The magnitude of the terms decreases until $n = 5$, after which the size of successive terms increases, causing the series to diverge. 

Typically, truncating the series at the value of $n$ that minimises $\eps^{2n}|u_n|$ produces the most accurate approximation that can be achieved by the series. This value of $n$ is often known as the ``optimal truncation point'', which we denote as $N_{\mathrm{opt}}$. 

Divergence of an asymptotic series, such as that seen in Figure \ref{Fig_seriesterms},  indicates that the solution contains  exponentially small terms that cannot be described by the algebraic series. Behaviour that occurs on this scale is smaller than any algebraic term in the limit that $\epsilon \to 0$, and hence cannot be captured by the series expression. 

Exponential asymptotic methods are built on the observation that truncating a series optimally at $n = N_{\mathrm{opt}}$ produces a remainder that is exponentially small in the asymptotic limit \cite{Berry}. By rescaling to study the exact truncation remainder, we can directly calculate exponentially small components of the asymptotic expansion.

\subsection{Finding the Stokes Curves}\label{S.exactStokes}

The analysis presented here is a standard application of the exponential asymptotic method developed in \cite{Daalhuis} in order to study Stokes' phenomenon. Nonetheless, this section will provide a relatively complete summary of the procedure in order to make the analysis accessible.

The purpose of this analysis is to obtain the leading exponentially small correction term to the series \eqref{eq.series}, and to therefore reveal the effects of Stokes switching in the solution. We will concentrate on the contribution arising from the term explicitly shown in \eqref{eq.seriesterms}, and note that a complex conjugate contribution is also present in the solution, which we will include in the final expression.

We now truncate the power series \eqref{eq.series} after $N$ terms, giving
\begin{equation}\label{eq.tseries}
u(x; \epsilon) = \sum_{n=0}^{N-1} \epsilon^{2n}u_n(x) + R_N(x),
\end{equation}
where $R_N$ is the remainder. The optimal truncation point $N_{\mathrm{opt}}$
corresponds to minimising $\eps^{2N}|u_N|$. The most straightforward way to find $N_{\mathrm{opt}}$ is to determine the value at which consecutive terms have equal magnitude, or
\begin{equation}
\left|\frac{(4n)!(-1)^n \epsilon^{2n}}{2^{4n}(2n)!(x-\i)^{2n+{1/2}}} \right| = \left|\frac{(4n+4)!(-1)^{n+1}\epsilon^{2n+2}}{2^{4n+4}(2n+2)!(x-\i)^{2n+{5/2}}}\right|,
\end{equation}
Optimal truncation occurs after a large number of terms in the limit that $\eps \to 0$. We can make use of this to find that the optimal truncation point must satisfy $N_{\mathrm{opt}} \sim |x-\i|/2\epsilon$ as $\eps \to 0$. We therefore set 
\begin{equation}\label{e.Nopt}
N_{\mathrm{opt}} = \frac{|x-\i|}{2\epsilon} + \omega,
\end{equation}
where $\omega \in [0,1)$ is chosen such that $N_{\mathrm{opt}}$ is an integer. This expression for $N_{\mathrm{opt}}$ is consistent with Figure \ref{Fig_seriesterms} which has $|x - \i| = 1$ and $\epsilon = 0.1$, corresponding to $N_{\mathrm{opt}} = 5$. We will eventually use this value for $N$ in \eqref{eq.tseries}, but for algebraic simplicity we will not make the substitution immediately.

We substitute the truncated series \eqref{eq.tseries} into the governing equation \eqref{eq.ode} to obtain
\begin{equation}
\sum_{n=0}^{N-1} \epsilon^{2n+2} u_n'' + \sum_{n=0}^{N-1}\epsilon^{2n} u_n + \epsilon^2 R_N'' + R_N =  \frac{1}{\sqrt{x+\i}} + \frac{1}{\sqrt{x-\i}}.
\end{equation}
Using  \eqref{eq.seriesterms} and \eqref{eq.recur} to simplify gives
\begin{equation}
 \epsilon^2 R_N'' + R_N =  \epsilon^{2N}u_N = \frac{ (4N)!(-1)^N\epsilon^{2N}}{2^{4N}(2N)!(x-\i)^{2N+{1/2}}},\label{eq.RN1}
\end{equation}
We will solve this differential equation using variation of parameters. The first step is to determine the homogeneous solutions to \eqref{eq.RN1}, given by
\begin{equation}
R_N(x) = K_1 \e^{\i x/\epsilon} \qquad \mathrm{and} \qquad R_N(x) = K_2 \e^{-\i x/\epsilon},
\end{equation}
where $K_1$ and $K_2$ are arbitrary constants. The next step is to permit the arbitrary constants to vary in $x$, and then consider the full differential equation in \eqref{eq.RN1}. We will find that for \eqref{eq.RN1} the second solution is the relevant one, so we set
\begin{equation}\label{eq.RNS}
R_N(x) = K_2(x)\e^{-\i x/\epsilon} = \mathcal{S}(x)\e^{-\i (x-\i)/\epsilon},
\end{equation}
where $K_2(x) = \mathcal{S}(x)\e^{-1/\epsilon}$, with the scaling chosen so that the exponent in \eqref{eq.RNS} is equal to 0 when $x = \i$. When written in this form, $\mathcal{S}(x)$ is known as the Stokes multiplier. It is constant except within a narrow region surrounding a Stokes curve. In this region, $\mathcal{S}$ varies rapidly from zero on one side of the curve to a nonzero value on the other. This rapid variation generates Stokes switching in the solution.

Substituting \eqref{eq.RNS} into \eqref{eq.RN1} and simplifying gives
\begin{equation}\label{eq.RN1b}
-2 \i \epsilon \mathcal{S}' \e^{-\i(x-\i)/\epsilon} = \frac{ (4N)!(-1)^N\epsilon^{2N}}{2^{4N}(2N)!(x-\i)^{2N+{1/2}}} .
\end{equation}
Recall that the optimal truncation given in \eqref{e.Nopt} depends on $|x-\i|$. This suggests that the transformation $\i(x - \i) = r\e^{\i\theta}$ will simplify this expression in a useful way, giving $N_{\mathrm{opt}} = r/2\epsilon + \omega$. As in \cite{Daalhuis}, we then fix $r$ and consider only the faster variation that occurs in the angular direction $\theta$. From the chain rule,
\begin{equation}
\diff{}{x} = \frac{\e^{-\i\theta}}{r}\diff{}{\theta}.
\end{equation}
Applying this transformation to \eqref{eq.RN1} and rearranging gives
\begin{equation}\label{eq.dSdt}
 \diff{\mathcal{S}}{\theta} = -\frac{\i\sqrt{\i}r\e^{\i\theta}(4N)!\epsilon^{2N-1}}{2^{4N+1}(2N)!(r\e^{\i\theta})^{2N+1/2}}\exp\left(\frac{r\e^{\i\theta}}{\epsilon}\right) \qquad \mathrm{as} \qquad \epsilon \to 0.
\end{equation}
We now finally substitute in the optimal truncation value $N_{\mathrm{opt}} = r/2\epsilon + \omega$. Making use of Stirling's formula in the limit that $\eps \to 0$ gives
\begin{equation}
\diff{\mathcal{S}}{\theta} \sim -\frac{\i\sqrt{2\i r}}{\epsilon}\exp\left(\frac{r}{\epsilon}\left(\e^{\i\theta} - 1 - \i\theta\right) + \i\theta\left(\frac{1}{2} + \omega\right)\right)\qquad \mathrm{as} \qquad \eps \to 0.
\end{equation}
As promised, this expression is exponentially small, except in the neighbourhood of $\theta = 0$, where the exponential term is algebraic. Hence, $\mathcal{S}$ varies rapidly in this region, indicating that $\theta = 0$ is a Stokes curve. 

\begin{figure}[tb]
\centering
 \begin{tikzpicture}
 [x=30,>=stealth,y=30]
 
 \fill[blue!10] (0,-2.25) -- (2.25,-2.25) -- (2.25,2.25) -- (0,2.25) -- cycle;

\draw[thick, decoration = {zigzag,segment length = 2mm, amplitude = 0.75mm}, decorate] (0,1) -- (0,2.25);
\draw[thick, decoration = {zigzag,segment length = 2mm, amplitude = 0.75mm}, decorate] (0,-1) -- (0,-2.25);

\draw[->] (-2.5,0) -- (2.5,0) node[right] {\scriptsize{$\mathrm{Re}(x)$}};
\draw[->] (0,-2.5) -- (0,2.5) node[right] {\scriptsize{$\mathrm{Im}(x)$}};

\draw[line width=1.5pt,red] (0,1) -- (0,-1);

\draw[line width=2pt] (0.1,0.9) -- (-0.1,1.1);
\draw[line width=2pt] (0.1,1.1) -- (-0.1,0.9);
\draw[line width=2pt] (0.1,-0.9) -- (-0.1,-1.1);
\draw[line width=2pt] (0.1,-1.1) -- (-0.1,-0.9);
\node at (0,1) [xshift=-3pt, left] {\scriptsize{i}};
\node at (0,-1) [xshift=-3pt, left] {\scriptsize{$-$i}};

\draw (3, 1.5) -- (3, 1.8) -- (3.3, 1.8) -- (3.3,1.5) -- cycle;
 \fill[blue!10] (3, 0.8) -- (3, 1.1) -- (3.3, 1.1) -- (3.3,0.8) -- cycle;
\draw (3, 0.8) -- (3, 1.1) -- (3.3, 1.1) -- (3.3,0.8) -- cycle;

\node[right,xshift=3pt] at (3.3,1.65) {\scriptsize{$u_{\mathrm{exp}} = 0$}};
\node[right,xshift=3pt] at (3.3,0.95)  {\scriptsize{$u_{\mathrm{exp}} \sim \sqrt{\frac{4\pi}{\epsilon}}\mathrm{e}^{-1/\epsilon} \cos\left(\frac{x}{\epsilon} + \frac{\pi}{4}\right)$}};

\draw[red,line width=1.5pt] (3,-0.667) -- (3.3,-0.667) node[right,black] {\scriptsize{Stokes curves}};
\draw[thick, decoration = {zigzag,segment length = 2mm, amplitude = 0.75mm}, decorate] (3,-1.333) -- (3.3,-1.333) node[right,black] {\scriptsize{Branch cut}};
\draw[line width=2pt] (3.05,-1.9) -- (3.25,-2.1);
\draw[line width=2pt] (3.25,-1.9) -- (3.05,-2.1);
\node[right] at (3.2,-2) {\scriptsize{Branch point}};

\node at (-1.125,-0.5) {\scriptsize{No exp.}};
\node at (-1.125,-0.75) {\scriptsize{terms}};
\node at (1.125,-0.5) {\scriptsize{Exp. small}};
\node at (1.125,-0.75) {\scriptsize{corrections}};

\draw[gray] (2.25,2.25) -- (-2.5,2.25) -- (-2.5,-2.25) -- (2.25,-2.25) -- cycle;

 \end{tikzpicture}
\caption{Stokes' phenomenon in the exact solution to \eqref{eq.ode} with boundary conditions \eqref{eq.bc}. The leading-order solution contains branch points at $x = \pm \i$, with the branch cuts extending vertically away from the real axis. The branch points generate Stokes curves, which connect the two points. On the left-hand side of the Stokes curves, there are no exponential contributions. On the right-hand side of the Stokes curves, the solution contains exponentially small oscillations of the form given in \eqref{e.uexp}.}\label{Fig:truestokes}
\end{figure}
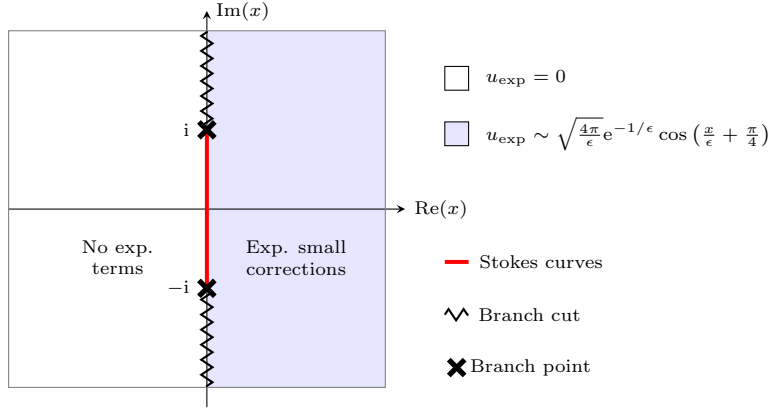

Recall that $\i(x-\i) = r\e^{\i\theta}$, so $\theta = 0$ corresponds to  $x - \i$ taking negative imaginary values, or a vertical line extending downwards from the branch point $x = \i$. If we perform a corresponding analysis for the complex conjugate term in \eqref{eq.seriesterms}, we identify that there is a second Stokes curve extending vertically upwards from the branch point $x = -\i$. This Stokes structure is illustrated in Figure \ref{Fig:truestokes}.

If we cross the Stokes line along $\mathrm{Re}(x) =0$, we expect that there will be an exponentially small jump in the asymptotic solution. 

\subsection{Exponentially Small Contribution}\label{S.exactMAE}

To determine the behaviour that appears as the Stokes curve is crossed, we define a new inner variable $\theta = \eps^{1/2}\vartheta$, which gives
\begin{equation}
\diff{\mathcal{S}}{\vartheta}  \sim -\i\sqrt{\frac{2\i r}{\epsilon}}\e^{-r\vartheta^2/2} \qquad \mathrm{as} \qquad \eps \to 0.
\end{equation}
By comparing $\vartheta$ with the original coordinates, we see that $x < 0$ corresponds to $\vartheta \to -\infty$, while $x > 0$ corresponds to $\vartheta \to \infty$. Hence, the jump across the Stokes curve is given by
\begin{equation}\label{eq.jump}
[\mathcal{S}]_-^+ = \lim_{\vartheta\to +\infty} \mathcal{S} -\lim_{\vartheta\to -\infty} \mathcal{S} \sim -\i\sqrt{\frac{2\i r}{\epsilon}}\int_{-\infty}^{\infty}\e^{-r\vartheta^2/2}\d \vartheta = -2 \i\sqrt{\frac{\pi\i}{\epsilon}}.
\end{equation}
If the value of $\mathcal{S}$ changes by \eqref{eq.jump} as the Stokes curve is crossed, the exponentially small contribution to the solution for $x > 0$ is given by
\begin{equation}\label{eq.RNjump}
[R_N]_-^+ \sim -2 \i\sqrt{\frac{\pi\i}{\epsilon}}\e^{-\i(x-\i)/\epsilon}.
\end{equation}
If we add in the complex conjugate term, we obtain the jump in the exponentially small terms $u_{\mathrm{exp}}$ as the Stokes curve is crossed from left to right,
\begin{equation}\label{e.uexpjump}
[u_{\mathrm{exp}}]_-^+ \sim  -2 \i\sqrt{\frac{\pi\i}{\epsilon}}\,\e^{-\i(x-\i)/\epsilon} + \mathrm{c.c.} = 4\sqrt{\frac{\pi}{\epsilon}}\,\e^{-1/\epsilon}\cos\left(\frac{x}{\epsilon}+\frac{\pi}{4}\right), \qquad \mathrm{as} \qquad \eps \to 0.
\end{equation}
Note that this expression has constant amplitude. By comparing this with the boundary conditions in \eqref{eq.bc}, we see that there cannot be any finite-amplitude oscillations present as $x \to -\infty$. This indicates that $u_{\mathrm{exp}} = 0$ on the left-hand side of the Stokes curve ($x < 0$), and that the exponentially small behaviour appears as the Stokes curve is crossed from left to right.
On the right-hand side of the Stokes curve, the exponentially small solution contribution $u_{\mathrm{exp}}$ must be
\begin{equation}\label{e.uexp}
u_{\mathrm{exp}} \sim  4\sqrt{\frac{\pi}{\epsilon}}\,\e^{-1/\epsilon}\cos\left(\frac{x}{\epsilon}+\frac{\pi}{4}\right), \qquad x > 0.
\end{equation}
This behaviour is illustrated in Figure \ref{Fig:truestokes}. 

\section{AAA Approximation}

The analysis in Section \ref{S.expasymp} was a relatively straightforward application of the method of \cite{Daalhuis}. It centered on the observation that singularities of the leading-order solution in the complex plane, such as the branch points at $x = \pm\i$ in \eqref{eq.LO}, generate Stokes curves. We can then optimally truncate and rescale the equation to study the exponentially small contributions in the neighbourhood of these terms.

\subsection{AAA Approximation for $u_0$}\label{S:AAA_u0}

In order to simulate a problem where we only possess a numerical description of the leading-order solution $u_0$ from \eqref{eq.LO}, we first define a set of points $x_j$ separated by some $\Delta x$. We then evaluate \eqref{eq.LO} at each point to obtain $u_0(x_j)$. We will use this set of points $x_j$ and function values $u_0(x_j)$ as the basis for a AAA rational approximation.

We apply the AAA algorithm to obtain a rational approximation $\hat{u}_0(x)$, given by
\begin{equation}\label{e.LOapp}
u_0 \approx \hat{u}_0 = \sum_{r=0}^{m} \frac{a_r}{x - p_r}.
\end{equation}
For example, if we choose $x_j$ to be evenly distributed in the interval $[-4,4]$ with $\Delta x = 0.1$, applying the AAA algorithm with an error tolerance of $10^{-12}$ gives a rational function with 15 poles. The poles and residues are are shown in Table \ref{PoleTable}. We have given each pair of poles a designation, so that we may refer to them later. 

\begin{table}
\begin{center}
\begin{tabular}{ |cccc| } 
 \hline
 Pole Index: $r$ & Pole Pair Name & Pole locations: $p_r$ & Pole residues: $a_r$  \\ 
 \hline
  
 \rowcolor{gray!15} 1, 2&Pair 1& $-0.0015 \pm 1.0178\i$ & $0.1256 \mp 0.1155\i$  \\ 
  
  3, 4&Pair 2&$-0.0142 \pm 1.1647\i$ & $0.1319 \mp 0.1204\i$ \\
  
  \rowcolor{gray!15} 5, 6&Pair 3&$-0.0444 \pm 1.4872\i$ & $0.1458 \mp 0.1312\i$ \\
  
  7, 8&Pair 4&$-0.1052 \pm 2.0575\i$ & $0.1725 \mp 0.1513\i$ \\
  
 \rowcolor{gray!15} 9, 10&Pair 5&$-0.2369 \pm 3.0523\i$ & $0.2283 \mp 0.1898\i$ \\
  
   11, 12&Pair 6&$-0.6072 \pm 5.0113\i$ & $0.3780 \mp 0.2791\i$ \\
  
 \rowcolor{gray!15}  13, 14&Pair 7&$-2.5372 \pm10.7492\i$ & $1.1038 \mp 0.6053\i$ \\
   
 15 &Unpaired & $-6.6066 + 0.0000\i$ & $-0.0003 + 0.0000\i$ \\
 \hline
\end{tabular}
\end{center}
\caption{Poles and residues in the AAA approximation of $u_0$, which contains branch points at $x = \pm\i$. The approximation $\hat{u}_0$ was produced by sampling $u_0$ on the domain $x \in [-4,4]$ at intervals of $\Delta x = 0.1$ with an error tolerance of $10^{-12}$. The poles in $\hat{u}_0$ accumulate at $x = \pm\i$, or the true branch points of $u_0$. An unpaired pole appears on the real axis outside the sampling interval, which we discard as a numerical artefact. The remaining poles occur in complex conjugate pairs. The pole $p_1$ is the nearest pole to the true branch point at $x = \i$, and this pole will play an important role in subsequent analysis.}\label{PoleTable}
\end{table}

Note that there is a pole located on the real axis at $x \approx -6.066$. The AAA algorithm sometimes generates poles on the real axis that lie outside of the approximation interval. These poles have no practical effect on the solution, and we will omit their contributions from the sum \eqref{e.LOapp}. The remaining poles appear as complex conjugate pairs, as the solution is real when $x$ is real. 

The pole pairs are not restricted to the imaginary axis, which is where we defined the branch cuts to be in the true solution $u_0$ \eqref{eq.LO}. In fact, the AAA algorithm places the poles along a curve, representing the branch cut, which we are not free to prescribe. This is a common feature of rational approximation methods, and was studied in depth for Pad\'{e} approximation in \cite{stahl1997convergence}, using a quantity known as the ``condenser capacity'' in the approximation. More recently, an explanation of the branch curve selection for AAA that built on these ideas was presented in \cite{trefethen2023numerical}.

It is possible to adjust the rational approximation algorithm, for example by using the \texttt{minimax} algorithm from the chebfun package, to obtain an approximation which forces the poles to align along the imaginary axis. 

\begin{figure}[tb]
    \centering
    \subfloat[Re($u_0$)]{
    \includegraphics[width=0.4\textwidth]{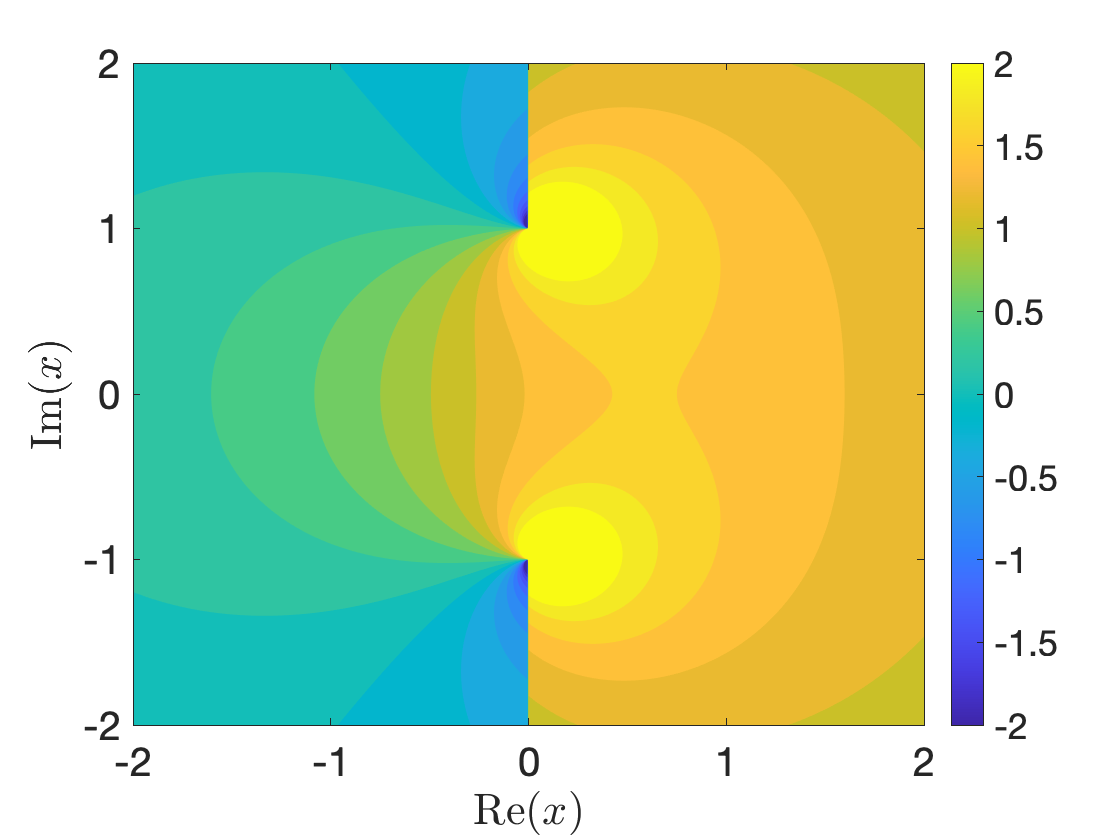}
    }
    \subfloat[Im($u_0$)]{
    \includegraphics[width=0.4\textwidth]{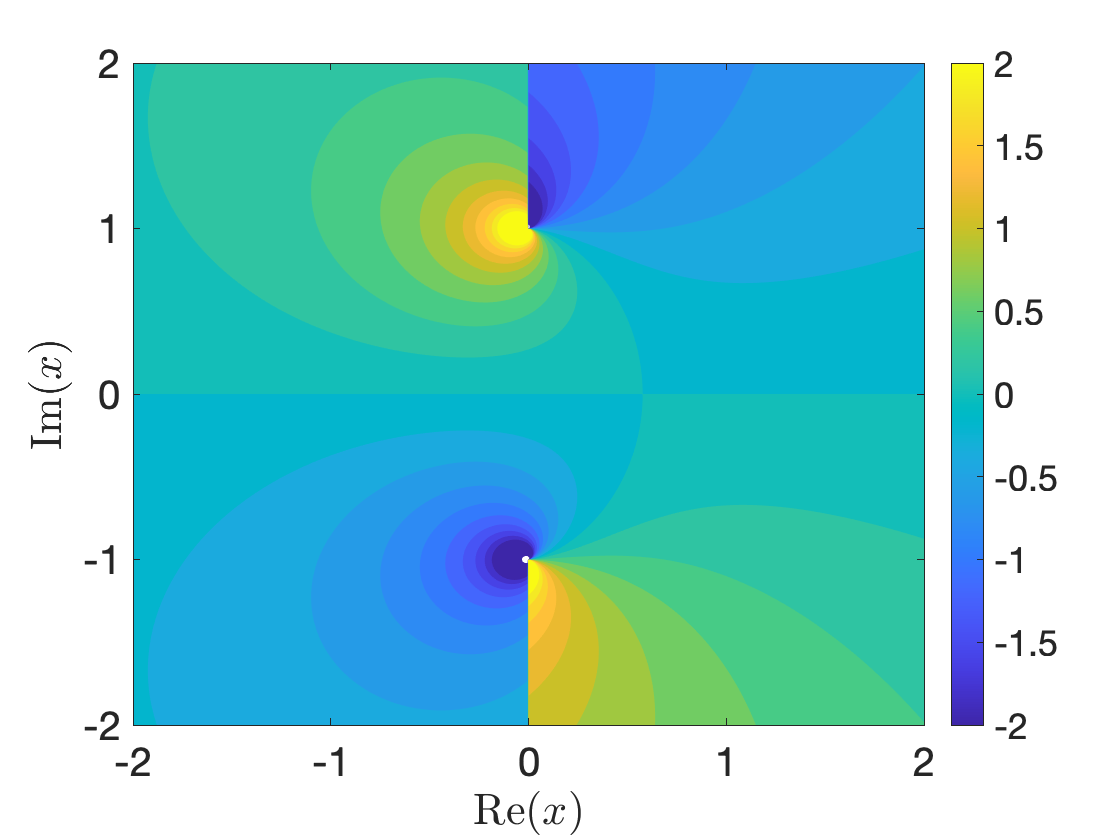}
    }

    \subfloat[Re($\hat{u}_0$)]{
    \includegraphics[width=0.4\textwidth]{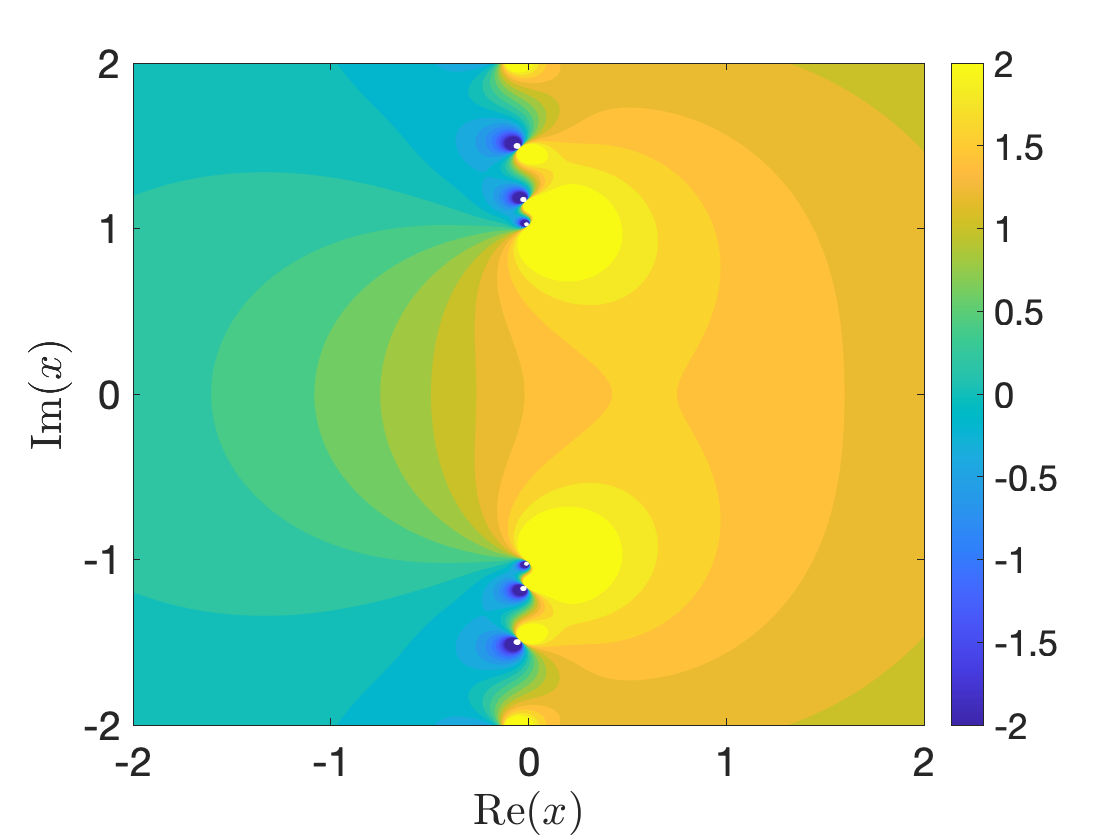}
    }
    \subfloat[Im($\hat{u}_0$)]{
    \includegraphics[width=0.4\textwidth]{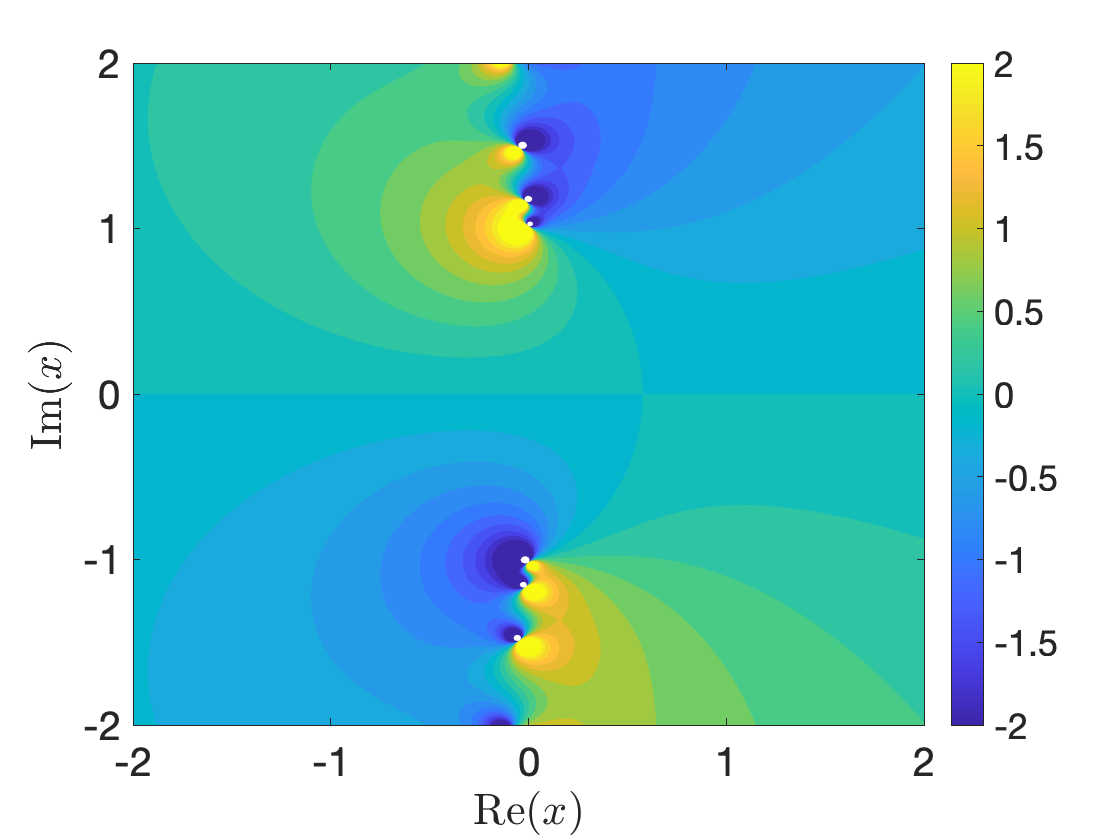}
    }
    \caption{The real and imaginary parts of the true leading-order solution $u_0$ and the approximated leading-order solution $\hat{u}_0$ described in Table \ref{PoleTable}. The two expressions are visually indistinguishable except on the imaginary axis. The function $u_0$ contains vertical branch cuts. The function $\hat{u}_0$ is a rational function which can only contain simple poles; these poles are arranged in such a way that they approximate the effect of a branch cut in the solution.}
    \label{fig:complexplots}
\end{figure}

In Figure \ref{fig:complexplots}, we present the real and imaginary parts of the true leading-order solution $u_0$, and the approximated leading-order solution $\hat{u}_0$. Note that the true solution possesses vertical branch cuts originating at $\pm\i$, while the approximated solution possesses simple poles that simulate the effect of a branch cut. The solutions appear visually identical except in a region surrounding the branch cut/poles. In this, and all subsequent analysis, we denote $p_1$ as the pole nearest to $x = \i$.

\begin{figure}[tb]
    \centering
    \subfloat[$\log_{10}|u_0 - \hat{u}_0|$]{
    \includegraphics[width=0.65\textwidth]{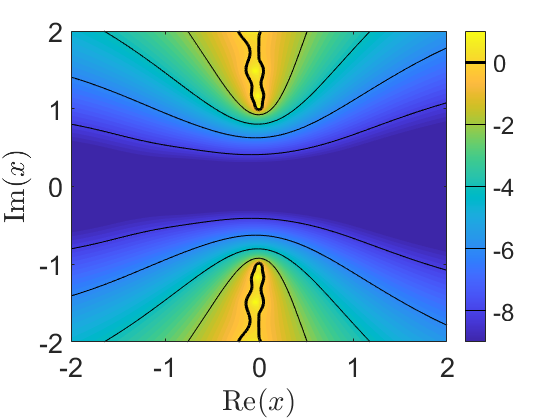}
    }
    \caption{This figure shows the error $|u_0 - \hat{u}_0|$, using $u_0$ and $\hat{u}_0$ from Figure \ref{fig:complexplots}. Note that the error is extremely small except in a region near the imaginary axis, where the true branch cut lies.}
    \label{fig:complexerror}
\end{figure}

In Figure \ref{fig:complexerror}, we present the error between the two on a logarithmic plot. This confirms that the approximation error is small except in a region surrounding the branch cut/poles, where the error becomes significant. Given that the approximation is highly accurate everywhere on the sample domain except near the branch cut, we hope that the exponential asymptotic analysis will produce equivalently accurate results.

\subsection{Asymptotic Power Series}

In effect, we are using the AAA approximation in place of the inhomogeneous term in \eqref{eq.ode}. Hence, we are determining the asymptotic behaviour of
\begin{equation}\label{eq.AAAode}
\epsilon^2 \hat{u}''(x) + \hat{u}(x) = \sum_{r=0}^{m} \frac{a_r}{x - p_r}.
\end{equation}
By expanding the solution as a series \eqref{eq.series} and matching powers of $\eps$, we can again obtain a recurrence relation for the series terms $u_n$. The recurrence relation gives
\begin{align}\label{e.AAAu0}
    \hat{u}_0(x) = \sum_{r=0}^{m} \frac{a_r}{x - p_r}, \qquad \hat{u}_n(x) = -\hat{u}_{n-1}''(x).
\end{align}
Solving this recurrence relation gives
\begin{equation}\label{e.uhatn}
    \hat{u}_n = \sum_{r=0}^{m} \frac{a_r (2n)!}{(x - p_r)^{2n+1}}.
\end{equation}
Each of the poles $p_r$ will generate a Stokes curve, and lead to an exponentially small contribution appearing in the solution.

\subsection{Stokes Curves and Exponentially Small Terms}

To determine the location of the Stokes curves, and the quantity that appears as they are crossed, we proceed with an almost identical analysis to that of Sections \ref{S.exactStokes} and \ref{S.exactMAE}, with the only significant difference being the form of the series terms in \eqref{e.uhatn}. 

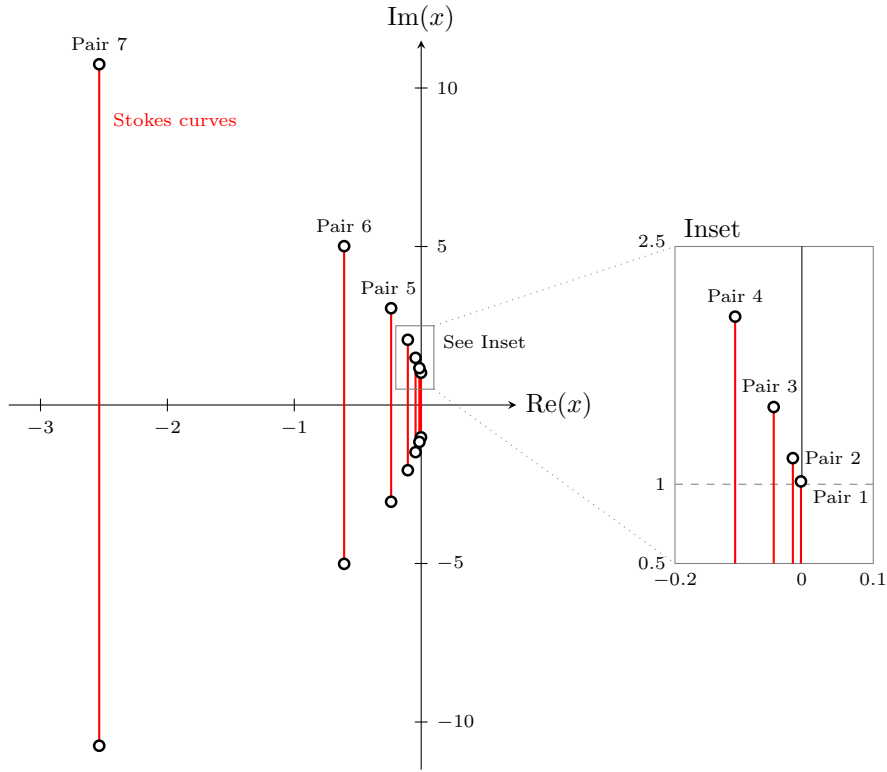
\begin{figure}[tb]
\centering
 \begin{tikzpicture}
 [x=12,>=stealth,y=12]

\draw[->] (-3.25*4,0) -- (0.75*4,0) node[right] {$\mathrm{Re}(x)$};
\draw[->] (0,-11.5) -- (0,11.5) node[above] {$\mathrm{Im}(x)$};

\draw[thick,red] (-2.5372*4,10.7492) --(-2.5372*4,-10.7492);
\draw[thick,red]  (-0.6072*4,5.0113) --(-0.6072*4,-5.0113) ;
\draw[thick,red]  (-0.2369*4,3.0523) --(-0.2369*4,-3.0523);
\draw[thick,red]  (-0.1052*4,2.0575) --(-0.1052*4,-2.0575) ;
\draw[thick,red]  (-0.0444*4,1.4872) --(-0.0444*4,-1.4872);
\draw[thick,red]  (-0.0142*4,1.1647) --(-0.0142*4,-1.1647) ;
\draw[thick,red]  (-0.0015*4,1.0178) --(-0.0015*4,-1.0178);

\fill (-0.0015*4,1.0178) circle (0.2); %node[left] {\scriptsize{Pair 1}};
\fill[white] (-0.0015*4,1.0178) circle (0.12); %node[left] {\scriptsize{Pair 1}};
\fill (-0.0142*4,1.1647) circle (0.2); %node[left] {\scriptsize{Pair 2}};
\fill[white] (-0.0142*4,1.1647) circle (0.12); %node[left] {\scriptsize{Pair 2}};
\fill (-0.0444*4,1.4872) circle (0.2); %node[left] {\scriptsize{Pair 3}};
\fill[white] (-0.0444*4,1.4872) circle (0.12); %node[left] {\scriptsize{Pair 3}};
\fill (-0.1052*4,2.0575) circle (0.2); %node[left] {\scriptsize{Pair 4}};
\fill[white] (-0.1052*4,2.0575) circle (0.12); %node[left] {\scriptsize{Pair 4}};
\fill (-2.5372*4,10.7492) circle (0.2) node[above,yshift=2pt] {\scriptsize{Pair 7}};
\fill[white] (-2.5372*4,10.7492) circle (0.12);% node[above] {\scriptsize{Pair 7}};
\fill (-0.6072*4,5.0113) circle (0.2) node[above,yshift=2pt] {\scriptsize{Pair 6}};
\fill[white] (-0.6072*4,5.0113) circle (0.12);% node[above] {\scriptsize{Pair 6}};
\fill (-0.2369*4,3.0523) circle (0.2) node[above,yshift=2pt,xshift=-1pt] {\scriptsize{Pair 5}};
\fill[white] (-0.2369*4,3.0523) circle (0.12);% node[above] {\scriptsize{Pair 5}};

\fill (-0.0015*4,-1.0178) circle (0.2) ;
\fill[white]  (-0.0015*4,-1.0178) circle (0.12) ;
\fill (-0.0142*4,-1.1647) circle (0.2);
\fill[white]  (-0.0142*4,-1.1647) circle (0.12);
\fill (-0.0444*4,-1.4872) circle (0.2) ;
\fill[white]  (-0.0444*4,-1.4872) circle (0.12) ;
\fill (-0.1052*4,-2.0575) circle (0.2);
\fill[white]  (-0.1052*4,-2.0575) circle (0.12);
\fill (-2.5372*4,-10.7492) circle (0.2);
\fill[white] (-2.5372*4,-10.7492) circle (0.12);
\fill (-0.6072*4,-5.0113) circle (0.2);
\fill[white]  (-0.6072*4,-5.0113) circle (0.12);
\fill (-0.2369*4,-3.0523) circle (0.2);
\fill[white]  (-0.2369*4,-3.0523) circle (0.12);

\draw[gray]  (8,5) node[left,yshift=2pt,black] {\scriptsize{$2.5$}}   node[above right,black] {{Inset}} -- (8,-5) node[below,black] {\scriptsize{$-0.2$}} node[left,black] {\scriptsize{$0.5$}}  -- (14.25,-5) node[below,black] {\scriptsize{$0.1$}} -- (14.25,5) -- cycle;
\draw[gray] (-0.2*4,0.5) -- (0.1*4,0.5)  -- (0.1*4,2.5) node[below right,black] {\scriptsize{See Inset}} -- (-0.2*4,2.5) -- cycle;
\draw (12,-5+5*1.0178-5*0.5) -- (12,5);
\draw[gray,dashed] (14.25,-2.5) -- (8,-2.5) node[black,left] {\scriptsize{$1$}};
\node[below] at (12,-5) {\scriptsize{$0$}}; 

\draw[gray,dotted] (0.1*4,2.5) -- (8,5);
\draw[gray,dotted] (0.1*4,0.5) -- (8,-5);

\draw[thick,red] (8-0.0015*20+0.2*20,-5+5*1.0178-5*0.5) -- (8-0.0015*20+0.2*20,-5);
\fill (8-0.0015*20+0.2*20,-5+5*1.0178-5*0.5) circle (0.2) node[below right,xshift=1pt] {\scriptsize{Pair 1}};
\fill[white]  (8-0.0015*20+0.2*20,-5+5*1.0178-5*0.5) circle (0.12) ;

\draw[thick,red] (8-0.0142*20+0.2*20,-5+5*1.1647-5*0.5) -- (8-0.0142*20+0.2*20,-5);
\fill (8-0.0142*20+0.2*20,-5+5*1.1647-5*0.5) circle (0.2) node[right,xshift=1pt] {\scriptsize{Pair 2}};
\fill[white]  (8-0.0142*20+0.2*20,-5+5*1.1647-5*0.5) circle (0.12) ;

\draw[thick,red] (8-0.0444*20+0.2*20,-5+5*1.4872-5*0.5) -- (8-0.0444*20+0.2*20,-5);
\fill (8-0.0444*20+0.2*20,-5+5*1.4872-5*0.5) circle (0.2) node[above,yshift=2pt,xshift=-1.5pt] {\scriptsize{Pair 3}};
\fill[white]  (8-0.0444*20+0.2*20,-5+5*1.4872-5*0.5) circle (0.12) ;

\draw[thick,red] (8-0.1052*20+0.2*20,-5+5*2.0575-5*0.5) -- (8-0.1052*20+0.2*20,-5);
\fill (8-0.1052*20+0.2*20,-5+5*2.0575-5*0.5) circle (0.2) node[above,yshift=2pt] {\scriptsize{Pair 4}};
\fill[white]  (8-0.1052*20+0.2*20,-5+5*2.0575-5*0.5) circle (0.12) ;

\draw (-1*4,-0.2) node[below] {\scriptsize{$-1$}} -- (-1*4,0.2);
\draw (-2*4,-0.2) node[below] {\scriptsize{$-2$}} -- (-2*4,0.2);
\draw (-3*4,-0.2) node[below] {\scriptsize{$-3$}} -- (-3*4,0.2);
\draw (-0.2,5) -- (0.2,5) node[right] {\scriptsize{$5$}};
\draw (-0.2,10) -- (0.2,10) node[right] {\scriptsize{$10$}};
\draw (-0.2,-5) -- (0.2,-5) node[right] {\scriptsize{$-5$}};
\draw (-0.2,-10) -- (0.2,-10) node[right] {\scriptsize{$-10$}};

\node[red,right] at (-10,9) {\scriptsize{Stokes curves}};

 \end{tikzpicture}
\caption{Stokes structure in the solution to \eqref{eq.AAAode}, which uses $\hat{u}_0$ as the leading-order solution. The solution contains 7 pairs of poles, with locations given in Table \ref{PoleTable}. Each pair of poles generates Stokes curves that extend vertically from the poles, intersecting the real axis. As each Stokes curve is crossed from left to right, an exponentially small asymptotic contribution appears in the solution. Note that the poles accumulate near the true branch points of $u_0$ at $x =\pm\i$. We will later find that the largest exponential contributions arise from the poles that are nearest to $x = \pm\i$.}\label{Fig.AAAStokes}
\end{figure}

This analysis reveals that there are exponentially small contributions in the solution associated with each pole. These contributions appear across Stokes curves that extends vertically from the corresponding pole and intersect the real axis. This behaviour is shown for the example problem in Figure \ref{Fig.AAAStokes}. We will later find that the most significant exponential contributions in this example appear across the Stokes curves generated by pole pairs 1, 2, and 3.

Using an essentially identical set of steps to Section \ref{S.exactMAE}, we can determine the form of the exponentially small remainder. We will determine the contribution that appears across the Stokes curve generated by a pole in the upper-half plane located at at $x = p_r$. The Stokes curve extends vertically downwards from the pole along $\mathrm{Re}(x) = \mathrm{Re}(p_r)$ and intersects the real axis at $x = \mathrm{Re}(p_r)$. 

We denote the optimally-truncated remainder as $\hat{R}_N$. On the left-hand side of the Stokes curve, we have $\hat{R}_N = 0$. On the right-hand side of the Stokes curve, we find that
\begin{equation}\label{eq.RNapp}
\hat{R}_N \sim \frac{2 \pi a_r}{\epsilon}\,\mathrm{e}^{-\i(x-p_r)/\eps} \quad \mathrm{as} \quad \epsilon \to 0.
\end{equation}
Once all of the Stokes curves have been crossed from left to right, we find that the combined exponentially small contribution, which we denote as $\hat{u}_{\mathrm{exp}}$, is therefore given by
\begin{equation}\label{e.AAAuexp}
\hat{u}_{\mathrm{exp}} \sim \sum_{r=0}^{m}\frac{2 \pi a_r}{\epsilon}\,\mathrm{e}^{-\i(x-p_r)/\eps} \quad \mathrm{as} \quad \epsilon \to 0.
\end{equation}
For example, for the set of poles presented in Table \ref{PoleTable}, we find that all of the exponentially small contributions are present in the solution on the right-hand side of the Stokes curve that intersects the real line at $x = -0.0015$.

Note that the exponential contributions in \eqref{e.AAAuexp} associated with upper-half plane poles at $x = p_r$ decay exponentially as $\mathrm{Im}(p_r)$ increases (with the complex conjugate contributions exhibiting corresponding exponential decay). This suggests that the poles nearest to the real axis will produce the largest exponential contributions. 

\section{Comparison}

We have now calculated the exponentially small contributions that appear in the asymptotic solution to the differential equation \eqref{eq.ode} with boundary condition \eqref{eq.bc}, given by $u_{\mathrm{exp}}$ in \eqref{e.uexp}. We have also calculated the exponentially small contributions that appear in the asymptotic solution to the approximate problem \eqref{eq.AAAode} with the same boundary condition, given by $\hat{u}_{\mathrm{exp}}$ in \eqref{e.AAAuexp}. The remaining question is whether the approximate exponential contributions are able to accurately approximate the true exponential contributions.

\begin{figure}[tb]
\centering
 \begin{tikzpicture}
 [x=75,>=stealth,y=850]

\draw (0,-0.065) -- (5,-0.065) -- (5,0.065) -- (0,0.065) -- cycle;
\draw[dotted] (0,-0.05) node[left] {\scriptsize{$-0.10$}} -- (5,-0.05);
\draw[dotted] (0,-0.025) node[left] {\scriptsize{$-0.05$}}  -- (5,-0.025);
\draw[dotted] (0,0.05) node[left] {\scriptsize{$0.10$}}  -- (5,0.05);
\draw[dotted] (0,0.025) node[left] {\scriptsize{$0.05$}}  -- (5,0.025);
 \node at (0,0) [left] {\scriptsize{$0$}} ;
 \node at (0,-0.065) [below] {\scriptsize{$0$}} ;
 \node at (5,-0.065) [below] {\scriptsize{$5$}} ;
 \draw (1, -0.062) -- (1, -0.065) node[below] {\scriptsize{$1$}};
 \draw (2, -0.062) -- (2,-0.065) node[below] {\scriptsize{$2$}};
 \draw (3,  -0.062) -- (3,-0.065) node[below] {\scriptsize{$3$}};
 \draw (4,  -0.062) -- (4,-0.065) node[below] {\scriptsize{$4$}};

 \draw[->] (5,-0.065) -- (5.1,-0.065) node[right] {\scriptsize{$x$}};
 \draw[->] (0,0.065) -- (0,0.074) node[right] {\scriptsize{$u_{\mathrm{exp}}$}};
 
 \draw[blue] plot[] file {Data_pole1.txt};
 \draw[blue] plot[] file {Data_pole2.txt};
 \draw[blue] plot[] file {Data_pole3.txt};
 \draw[blue] plot[] file {Data_pole4.txt} node[right,yshift=-3pt] {\scriptsize{Pairs 4--7}};
 \draw[blue] plot[] file {Data_pole5.txt} node[right,yshift=3pt] {\scriptsize{Pair 3}};
 \draw[blue] plot[] file {Data_pole6.txt} node[right,yshift=2pt] {\scriptsize{Pair 2}};
\draw[blue] plot[] file {Data_pole7.txt} node[right,yshift=2pt] {\scriptsize{Pair 1}};

\draw[red, line width=1pt] plot[] file {Data_poles.txt} node[below right,yshift=5pt] {\scriptsize{$\hat{u}_{\mathrm{exp}}$}};
\draw[black, line width=1.75pt,dashed] plot[] file {Data_true.txt} node[above right] {\scriptsize{$u_{\mathrm{exp}}$}};

 \end{tikzpicture}
\caption{Exponentially small oscillations in the asymptotic solution to \eqref{eq.ode} and \eqref{eq.AAAode} for $\epsilon = 0.2$. The true exponential contribution $u_{\mathrm{exp}}$ is shown as a black dashed curve. The approximate exponential contribution $\hat{u}_{\mathrm{exp}}$ is shown as a red curve. This contribution was generated using the poles and residues from Table \ref{EpsTable}. The two curves are visually indistinguishable. The contribution $\hat{u}_{\mathrm{exp}}$ is the sum of contributions from each of the seven pole pairs. These contributions are shown individually as blue curves; it is apparent that the largest contributions arise from the pole pairs that are nearest to $x = \pm\i$ (ie. pole pairs 1, 2, and 3), with the amplitude of the contributions decaying as the distance of the pair from $x = \pm\i$ increases.}\label{Fig:e0p2}
\end{figure}
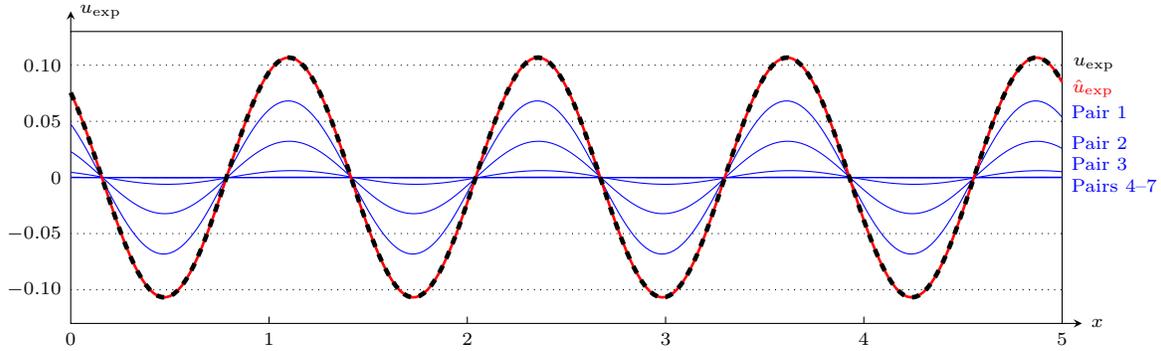

Note that \eqref{e.uexp} and \eqref{e.AAAuexp} cannot be identical as $\eps \to 0$, due to the following differences:
\begin{itemize}
    \item The algebraic prefactor of $u_{\mathrm{exp}}$ is proportional to $\eps^{-1/2}$, while the algebraic prefactor of the approximate exponential behaviour $\hat{u}_{\mathrm{exp}}$ is proportional to $\eps^{-1}$. This difference in exponent will be generic: the AAA rational function must have simple poles that generate algebraic prefactors that are proportional to $\eps^{-1}$, while more general singular points lead to prefactors whose algebraic power depends on the strength of the singularity.
    \item The exponential scaling of each expression is determined by the location of singular points in the leading-order solutions $u_0$ \eqref{eq.LO} and $\hat{u}_0$ \eqref{e.AAAu0}. The branch points in $u_0$ that produce \eqref{e.uexp} are located at $x = \pm\i$, while the poles in $\hat{u}_0$ that produce \eqref{e.AAAuexp} are located at $x = p_r$, where $|\mathrm{Im}(p_r)| > 1$. Hence the two expressions have different exponential decay as $\eps \to 0$, which can be seen by directly comparing the two expressions. 
\end{itemize}
These differences appear to suggest that the two expressions $u_{\mathrm{exp}}$ and $\hat{u}_{\mathrm{exp}}$ cannot possibly describe the same behaviour, due to the difference in the strength and position of the singular points between the exact and approximate leading-order behaviour. Indeed, the two expressions must be different in the limit that $\epsilon \to 0$, as all of the exponential terms in \eqref{e.AAAuexp} are smaller than the exponential terms in \eqref{e.uexp} in the asymptotic limit (even though the algebraic prefactor is larger). 

Despite this, we find that $\hat{u}_{\mathrm{exp}}$ is able to accurately approximate $u_{\mathrm{exp}}$ for some range of $\epsilon$. We compare the two contributions for our sample problem with $\epsilon = 0.2$ over a segment of the positive real axis in Figure \ref{Fig:e0p2}. The contributions $u_{\mathrm{exp}}$ and $\hat{u}_{\mathrm{exp}}$ appear indistinguishable in this figure.

In addition to the exponentially small contributions, Figure \ref{Fig:e0p2} also shows the contribution to $\hat{u}_{\mathrm{exp}}$ from each pole pair, and reveals that the largest contributions come from pole pairs 1--3. This indicates that pairs which are nearest to the branch points at $x = \pm\i$ also have the most significant effect on the exponentially-small behaviour of the solution, which is consistent with the exponential decay of the solution contributions as $|\mathrm{Im}(p_r)|$ increases.

\begin{figure}[tb]
\centering
 \begin{tikzpicture}
 [x=2500,>=stealth,y=100]

\draw[dotted] (0,1) node[left] {\scriptsize{$1$}} -- (0.15,1);
\draw (0.0016,0.5) -- (0,0.5) node[left] {\scriptsize{$0.5$}};
\draw[dotted] (0,0.5) -- (0.1,0.5);
\draw[dotted] (0.05,0) -- (0.05,1.15);
\draw[dotted] (0.1,0) -- (0.1,1.15);

\draw[blue, line width=1pt] plot[] file {Fig_t14_p10.txt};
\draw[red, line width=1pt] plot[] file {Fig_t13_p9.txt};
\draw[blue, line width=1.75pt,dotted] plot[] file {Fig_t12_p9.txt};
\draw[red, line width=1.75pt,dotted] plot[] file {Fig_t11_p8.txt};
\draw[black, line width=1pt] plot[] file {Fig_t10_p7.txt};
\draw[magenta, line width=1pt] plot[] file {Fig_t9_p7.txt};
\draw[black, line width=1.75pt,dotted] plot[] file {Fig_t8_p6.txt};
\draw[magenta, line width=1.75pt,dotted] plot[] file {Fig_t7_p5.txt};
\draw[black,dashed, line width=1.5pt] plot[] file {Fig_t6_p4.txt};

\fill[white] (0,0) -- (0.15,0) -- (0.15,-0.03) -- (0,-0.03) -- cycle;
\draw (0,0) node[left] {\scriptsize{$0$}} -- (0.15,0) -- (0.15,1.15) -- (0,1.15) -- cycle;

\draw(0.05,0.04) -- (0.05,0) node[below] {\scriptsize{$0.05$}};
\draw(0.1,0.04) -- (0.1,0) node[below] {\scriptsize{$0.10$}};
\node at (0,0) [below] {\scriptsize{$0$}};
\node at (0.15,0) [below] {\scriptsize{$0.15$}};

\draw[blue, line width=1pt] (0.1035,0.9) -- (0.11,0.9) node[right,black] {\scriptsize{$\mathrm{Tol.} = 10^{-14}$, 10 pairs}};
\draw[red, line width=1pt] (0.1035,0.8) -- (0.11,0.8) node[right,black] {\scriptsize{$\mathrm{Tol.} = 10^{-13}$, 9 pairs}};
\draw[blue, line width=1.75pt,dotted] (0.1035,0.7) -- (0.11,0.7) node[right,black] {\scriptsize{$\mathrm{Tol.} = 10^{-12}$, 9 pairs}};
\draw[red, line width=1.75pt,dotted] (0.1035,0.6) -- (0.11,0.6) node[right,black] {\scriptsize{$\mathrm{Tol.} = 10^{-11}$, 8 pairs}};
\draw[black, line width=1pt] (0.1035,0.5) -- (0.11,0.5) node[right,black] {\scriptsize{$\mathrm{Tol.} = 10^{-10}$, 7 pairs}};
\draw[magenta, line width=1pt] (0.1035,0.4) -- (0.11,0.4) node[right,black] {\scriptsize{$\mathrm{Tol.} = 10^{-9}$, 7 pairs}};
\draw[black, line width=1.75pt,dotted] (0.1035,0.3) -- (0.11,0.3) node[right,black] {\scriptsize{$\mathrm{Tol.} = 10^{-8}$, 6 pairs}};
\draw[magenta, line width=1.75pt,dotted] (0.1035,0.2) -- (0.11,0.2) node[right,black] {\scriptsize{$\mathrm{Tol.} = 10^{-7}$, 5 pairs}};
\draw[black,dashed, line width=1.5pt]  (0.11,0.1) node[right,black] {\scriptsize{$\mathrm{Tol.} = 10^{-8}$, 4 pairs}} -- (0.1035,0.1);

\draw[->] (0.15,0) -- (0.1535,0) node[right] {$\epsilon$};
\draw[->] (0,1.15) -- (0,1.15+0.0875) node[right] {Amplitude of $\hat{u}_{\mathrm{exp}}$/Amplitude of $u_{\mathrm{exp}}$};

\draw[black] plot[mark=*] (0.0450,1-0.1607);
\draw[magenta] plot[mark=*] (0.0327,1-0.1642);
\draw[black] plot[mark=*] (0.0264,1-0.1619);
\draw[magenta] plot[mark=*] (0.0193,1-0.1656);
\draw[black] plot[mark=*] (0.0178,1-0.1643);
\draw[red] plot[mark=*] (0.0151,1-0.1670);
\draw[blue] plot[mark=*] (0.0121,1-0.1671);
\draw[red] plot[mark=*] (0.0121,1-0.1671);
\draw[blue] plot[mark=*] (0.0096,1-0.1672);

\draw[blue] plot[mark=*] (0.055,0.1);
\draw[red] plot[mark=*] (0.058,0.1);
\draw[black] plot[mark=*] (0.061,0.1);
\draw[magenta] plot[mark=*] (0.064,0.1);
\node[right]  at (0.065,0.1) {\scriptsize{$\epsilon = |p_{\mathrm{n}} - \i|$}};

 \end{tikzpicture}
\caption{Ratio of the amplitudes of the approximate exponential contribution $\hat{u}_{\mathrm{exp}}$ and the true exponential contribution $u_{\mathrm{exp}}$ as $\epsilon$ is varied. If $\hat{u}_{\mathrm{exp}}$ is an accurate approximation of $u_{\mathrm{exp}}$, this ratio is close to 1. The two expressions have different behaviour as $\eps \to 0$, so it is impossible for the approximation to remain accurate indefinitely as $\eps$ is decreased. This is apparent in the figure as each approximation is accurate until some lower threshold value of $\eps$ is reached, beyond which the approximation becomes inaccurate. Increasing the tolerance, and hence the number of poles in the AAA approximation, reduces this threshold value of $\eps$. We also identify the value of $\epsilon$ equal to $|p_1 - \i|$ on each curve. This value of $\epsilon$ corresponds to a roughly constant value of the error in each curve, suggesting that the approximation error of the exponential terms is related the quantity $|p_1 - \i|$. }\label{Fig:AllEps}
\end{figure}
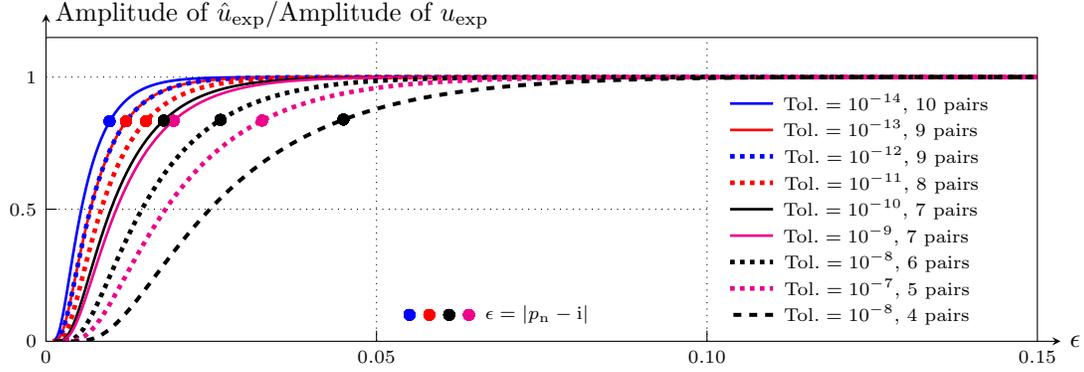

In Figure \ref{Fig:AllEps}, we compare the amplitude of $u_{\mathrm{exp}}$ and $\hat{u}_{\mathrm{exp}}$ over a range of $\epsilon$ values. We depict several different curves, each of which shows the results for a different AAA algorithm tolerance. 

For each curve, the ratio is approximately 1 (indicating that the approximation is accurate) until some threshold value of $\epsilon$ is reached from above; below this threshold value the ratio tends to 0. This is consistent with the algebraic expressions for the two terms, as the exponential contributions in $\hat{u}_{\mathrm{exp}}$ are smaller than the contributions in $u_{\mathrm{exp}}$ as $\eps \to 0$.

Reducing the error tolerance of the AAA algorithm has the effect of decreasing the value of $\epsilon$ for which the approximation of the exponential terms is accurate. Furthermore, this reduction does not happen in a uniform fashion. Changing the tolerance from $10^{-9}$ to $10^{-10}$ or from $10^{-12}$ to $10^{-13}$ does not change the number of poles in the approximate leading-order behaviour, and we see that this change has little effect on the threshold value of $\epsilon$.

It is apparent from Figure \ref{Fig:AllEps} that higher AAA tolerances produce approximations with more poles, that are accurate for smaller values of $\epsilon$. In Table \ref{EpsTable}, we show the effect of the error in branch cut prediction, or the difference in location between the true branch point (at $x = \i$) and the nearest pole in the approximation (at $x = p_1$), or $|p_1 - \i|$. We compare this with the relative error in the amplitude, measured by
\begin{equation}
    \textrm{Relative Error} = 1 - \frac{\textrm{Amplitude of } \hat{u}_{\mathrm{exp}}}{\textrm{Amplitude of } {u}_{\mathrm{exp}}}.
\end{equation}

\begin{table}
\begin{center}
\begin{tabular}{ |cccc| } 
 \hline
 AAA Tolerance & No. of Pairs & $\quad |p_1 - \i|\quad$ & Rel. Err. at $\eps = |p_1 - \i|$  \\ 
 \hline
 \rowcolor{gray!15} $10^{-6}$& 4&$0.0450$ & $0.1607$ \\
  $10^{-7}$& 5&$0.0327$ & $0.1642$ \\
  \rowcolor{gray!15}$10^{-8}$& 6&$0.0264$ & $0.1619$ \\
 $10^{-9}$&7&$0.0193$ & $0.1656$ \\
  \rowcolor{gray!15}$10^{-10}$&7&$0.0178$ & $0.1643$ \\
 $  10^{-11}$&8&$0.0151$ & $0.1670$ \\
 \rowcolor{gray!15}$10^{-12}$&9& $0.0121$ & $0.1670$ \\
 $ 10^{-13}$&9&$0.0121$ & $0.1670$ \\
 \rowcolor{gray!15} $10^{-14}$&10& $0.0096$ & $0.1671$ \\
   \hline
\end{tabular}
\end{center}
\caption{This table shows the distance between the true branch point in $u_0$ at $x = \i$ and the nearest pole in $\hat{u}_0$, denoted as $p_1$, for the different AAA approximations presented in Figure \ref{Fig:AllEps}. As the tolerance of the approximation increases, the distance $|p_1 - \i|$ decreases.  The relative approximation error is evaluated for $\eps = |p_1 - \i|$, and the values are roughly constant. This suggests that the approximation error depends on the quantity $|p_1 - \i|$, or how accurately the AAA approximation predicts the location of the true branch point.} \label{EpsTable}
\end{table}

Table \ref{EpsTable} suggests that, no matter what AAA tolerance is chosen, the relative error remains relatively similar at $\epsilon = |p_1 - \i|$. The points $\eps = |p_1 - \i|$ are also marked in Figure \ref{Fig:AllEps}. This figure shows that for $\eps \ll |p_1 - \i|$, the approximated amplitudes are inaccurate while for $\eps \gg |p_1 - \i|$, the approximate amplitudes are accurate. We conjecture that if there is a singularity in the true leading-order behaviour at $x = x_{\mathrm{s}}$ and the nearest pole to this point in the rational approximation is at $x = p_1$, the exponential asymptotic predictions will be accurate as long as $|p_1 - x_{\mathrm{s}}| \ll \epsilon$.

Obviously, this method is most valuable for problems in which $x_{\mathrm{s}}$ is not known beforehand, so $|p_1 - x_{\mathrm{s}}|$ is not known. It is possible to estimate $x_{\mathrm{s}}$ by carefully examining the accumulation rate of poles in the complex plane \cite{Trefethen2021}. We intend to study this in more detail in future work, and to explain the dependence of the threshold value on $|p_1 - x_{\mathrm{s}}|$.

\section{A Nonlinear Differential Equation}

From this analysis, it appears that rational approximation methods can be used to apply exponential asymptotic methods to linear ordinary differential equations, as long as care is taken to ensure that $\epsilon$ is not too small. The next natural question is whether we can extend these methods to nonlinear ordinary differential equations, using the exponential asymptotic method of \cite{Chapman}. 

This is a challenging problem, because we can no longer add the pole contributions in $\hat{u}_{\mathrm{exp}}$ independently. From \cite{trinh2015exponential}, we know that singularities in nonlinear problems that are near to each other (ie. within a neighbourhood whose width is proportional to a particular power of $\epsilon$) can interact, and that the resultant asymptotic behaviour changes due to these interactions. This is likely to occur for leading-order solutions generated using rational approximation methods due to the accumulation of poles near the true singular point. In future work, we will determine the asymptotic corrections to $u_{\mathrm{exp}}$ due to pole interactions in nonlinear problems.

In many problems it is possible to determine the strength of the singular points in $u_0$ using asymptotic balancing arguments, even if it is not possible to easily determine their location. In this case, another possible method for studying these problems is to apply a map to the sampled data so that the AAA algorithm is being fitted to a function with simple poles, then to invert the map so that it contains singularities with the correct order. 

Consider the following simple example:
%\jon{repetition}
\begin{equation}\label{e.nonlinode}
    \epsilon^2 u''(x) + u(x)^2 = \frac{81}{1024}\left(\frac{1}{\sqrt{x-\i}} + \frac{1}{\sqrt{x+\i}}\right)^2,
\end{equation}
with the boundary conditions in \eqref{eq.bc}. The constant $81/1024$ was chosen to make the exponential asymptotic analysis particularly convenient. The true leading-order behaviour is given by
\begin{equation}\label{e.nonlinAAAu0}
    u_0 = \frac{9}{32}\left(\frac{1}{\sqrt{x-\i}} + \frac{1}{\sqrt{x+\i}}\right).
\end{equation}
We can sample this leading-order solution on a set of points $x_j$ to obtain data points $u_0(x_j)$, as in Section \ref{S:AAA_u0}, and use this as the basis for an approximation $\hat{u}_0$. Using this, we can compute the exponential contributions that appear across each Stokes curve in the solution independently and add them together. If we do so, however, we obtain predictions of the exponentially small behaviour that are inaccurate, because the analysis does not take into account nonlinear interactions between pole contributions. 

We instead sample the leading-order solution on $x_j$ and then square the output, to obtain $u_0(x_j)^2$, which is equivalent to taking samples of
\begin{equation}\label{e.u02}
    u_0^2 = \frac{81}{1024}\left(\frac{1}{x-\i} + \frac{1}{x+\i} + \frac{2}{\sqrt{x^2+1}}\right).
\end{equation}
We can obtain a AAA rational approximation for $u_0^2$ based on this data, which we denote $\widehat{u_0^2}$. Finally, we can take the square root of this rational approximation to obtain 
\begin{equation}\label{e.nonlinAAALO}
    u_0 \approx \left(\widehat{u^2_0}\right)^{1/2} = \sqrt{\sum_{r=0}^m \frac{a_r}{x - p_r}}.
\end{equation}
This expression can then be used as the leading-order solution for an exponential asymptotic analysis using the methods from \cite{Chapman}. 
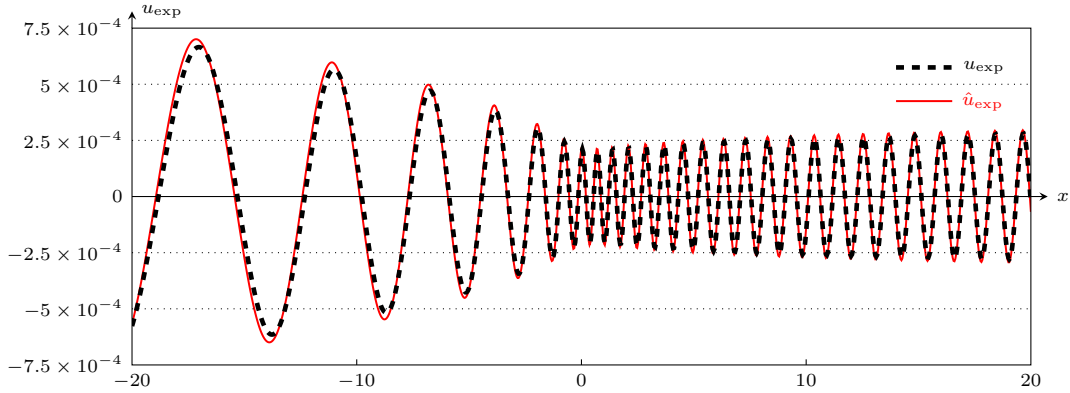
\begin{figure}[tb]
\centering
 \begin{tikzpicture}
 [x=8.5,>=stealth,y=8.5]
  \draw[red, line width=0.75pt] plot file {Nonlinear_R2.txt};
  \draw[black, line width=1.75pt,dashed] plot file {Nonlinear_R1.txt};

 \fill[white] (20,-7.5) -- (20,7.5) -- (20.3,7.5) -- (20.3,-7.5) -- cycle;
  \fill[white] (-20,-7.5) -- (-20,7.5) -- (-20.3,7.5) -- (-20.3,-7.5) -- cycle;
 
 \draw(-20,-7.5) -- (20,-7.5) -- (20,7.5) -- (-20,7.5) -- cycle;
 \draw[->] (-20,0) node[left] {\scriptsize{$0$}} -- (20.75,0) node[right] {\scriptsize{$x$}};
 \draw[->] (-20,7.5) -- (-20,8.25) node[right] {\scriptsize{$u_{\mathrm{exp}}$}};
 
 \draw[dotted] (-20,5) node[left] {\scriptsize{$5\times 10^{-4}$}} -- (20,5);
  \draw[dotted] (-20,-5) node[left] {\scriptsize{$-5\times 10^{-4}$}} -- (20,-5);
   \draw[dotted] (-20,2.5) node[left] {\scriptsize{$2.5\times 10^{-4}$}} -- (20,2.5);
  \draw[dotted] (-20,-2.5) node[left] {\scriptsize{$-2.5\times 10^{-4}$}} -- (20,-2.5);
  \node[left] at (-20,7.5) {\scriptsize{$7.5\times 10^{-4}$}};
  \node[left] at (-20,-7.5) {\scriptsize{$-7.5\times 10^{-4}$}};
  \node [below] at (-20,-7.5) {\scriptsize{$-20$}};
\node [below] at (20,-7.5) {\scriptsize{$20$}};
\draw (-10,-7.25) -- (-10,-7.5) node[below] {\scriptsize{$-10$}};
\draw (0,-7.25) -- (0,-7.5) node[below] {\scriptsize{$0$}};
\draw (10,-7.25) -- (10,-7.5) node[below] {\scriptsize{$10$}};

\draw[black, line width=1.75pt,dashed] (14,5.75) -- (16.5,5.75) node[right] {\scriptsize{$u_{\mathrm{exp}}$}};
\draw[red, line width=0.75pt] (16.5,4.25) node[right] {\scriptsize{$\hat{u}_{\mathrm{exp}}$}} -- (14,4.25) ;

 \end{tikzpicture}
\caption{Comparison of the true exponential terms and the AAA-approximated exponential terms for the nonlinear differential equation \eqref{e.nonlinode} with $\eps = 0.1$. The AAA-approximated exponential terms were obtained using the leading-order from \eqref{e.nonlinAAALO}. In order to show that this method captures nonlinear effects, we show the expression on the entire real axis, but note that this contribution is only actually present in the asymptotic solution for $x > 0$. The two contributions are visually similar but not completely identical, due to nonlinear effects caused by the subdominant branch points at $x = \pm\i$ in \eqref{e.u02}.}\label{Fig:Nonlin}
\end{figure}
In Figure \ref{Fig:Nonlin}, we present the true exponentially small terms $u_{\mathrm{exp}}$ and the approximated exponentially small terms $\hat{u}_{\mathrm{exp}}$ obtained using this method for $\epsilon = 0.1$. The rational approximation used was generated using sample points on $x\in[-10,10]$ with $\Delta x = 0.1$. We present the exponential contribution on $x \in [-20,20]$ to demonstrate that this method accurately captures nonlinear wave effects, but we actually expect this contribution to appear as the Stokes curve intersecting $x = 0$ is crossed from left to right, and therefore only be present in the asymptotic solution for $x > 0$. The qualitative match seen in this figure shows that the exponentially-small terms in \eqref{e.nonlinode} can be accurately approximated using this method.

The accuracy of this approximation compared to naively approximating the exponential terms based on $\hat{u}_0$ can be understood by comparing the poles and residues of the AAA approximation of $u_0$ and $u_0^2$, using the parameters from Figure \ref{Fig:Nonlin}. These are presented in Table \ref{NonlinTable}.

In $u_0$, the only true singularities are the branch points at $x = \pm\i$. In Table \ref{NonlinTable}, we see that the residue of each pole of $\hat{u}_0$ is similar in magnitude, as these poles approximate the effect of the branch cuts.  In $u_0^2$, there are simple poles at $x = \pm\i$, as well as branches that also originate at these points. The dominant behaviour for this expression is the simple poles.  In Table \ref{NonlinTable}, we see that the poles located closest to $\pm\i$ have a larger residue, with the other poles having a smaller residue. This is because the poles nearest to $\pm\i$ in the AAA approximation are reproducing the simple pole behaviour in the true expression for $u_0^2$, and the remaining poles are approximating the subdominant branch cut behaviour.

\begin{table}[tb]
\begin{center}
\begin{tabular}{ |p{1.5cm}cp{1.5cm}cc| } 
\hline
\multicolumn{1}{|c}{} & \multicolumn{2}{c}{$\hat{u}_0$}  & \multicolumn{2}{c|}{$\left(\widehat{u^2_0}\right)^{1/2}$}  \\ 
 \hline
 \rowcolor{gray!15}  &  $p_r$ &  \hspace{6pt}$|a_r|$ &  $p_r$ &  $|a_r|$  \\ 
 
 \rowcolor{white} Pair 1& $0.0020 \pm 1.0172\i$ & $0.2173$& $0.0018 \pm 1.0004\i$ & $0.2845$  \\ 
  
  \rowcolor{gray!15} Pair 2&$0.0180 \pm 1.1582\i$ & $0.2210$& $-0.0067 \pm 1.1204\i$ & $0.1342$ \\
  
  \rowcolor{white} Pair 3&$0.0507 \pm 1.4604\i$ & $0.2287$&$ -0.0207 \pm 1.3724\i$ & $0.1329$ \\
  
  \rowcolor{gray!15} Pair 4&$0.1014 \pm 1.9686\i$ & $0.2404$& $-0.0381 \pm 1.7959\i$ & $0.1335$ \\
 
 \rowcolor{white} Pair 5&$0.1728\pm 2.7620\i$ & $0.2566$& $-0.0550 \pm 2.4573\i$ & $0.1330$ \\
  
 \rowcolor{gray!15}  Pair 6&$0.2726 \pm 3.9795\i$ & $0.2793$& $-0.0665 \pm 3.4701\i$ & $0.1374$ \\
 \hline
\end{tabular}
\end{center}
\caption{Comparison of the six pole pairs nearest to $x = \pm\i$ contained in a rational approximation for $u_0$ from \eqref{e.nonlinAAAu0} with the six nearest pole pairs in the approximation for $u_0^2$ from \eqref{e.u02}. The poles in $\hat{u}_0$ have residues with similar magnitude, as the singularities in $u_0$ are branch points at $x = \pm\i$. In $(\widehat{u}_0^2)^{1/2}$, the poles nearest to $x = \pm\i$ have a larger residue than the remaining poles, as it approximating the contribution from the simple poles in \eqref{e.u02}. The remaining poles mimic the behaviour of the subdominant branch cut.} \label{NonlinTable}
\end{table}

Because the strongest singularities in $u_0^2$ from \eqref{e.u02} are simple poles, the AAA approximation can be used as the basis for an accurate exponential asymptotic analysis. However, this expression does still contain branch cuts originating at $x = \pm\i$. This produces a string of additional poles in the AAA approximation. Hence, any higher-order corrections to the exponential terms, such as those studied in \cite{shelton2023exponential,chapman2005exponential}, will be inaccurate unless pole interaction effects are taken into account. A more thorough analysis of nonlinear differential equations using rational approximation methods is beyond the scope of this article and will be the subject of future work. 

\section{Conclusions and Discussion}

In this study, we applied exponential asymptotic methods from \cite{Daalhuis} to study Stokes' phenomenon in a linear ordinary differential equation in the small-$\epsilon$ limit. The leading-order solution $u_0$ \eqref{eq.LO} contains two branch points, located at $x = \pm \i$. These branch points produce Stokes curves, and oscillating exponentially small terms appear in the solution as these Stokes curves are crossed at $x = 0$ on the real axis. We calculated the form of these exponentially small oscillations, $u_{\mathrm{exp}}$ \eqref{e.uexp}.

We then repeated this process using a rational approximation for the leading-order behaviour. Instead of using the exact leading-order solution, we sampled the leading-order behaviour on a discrete set of points and used the AAA algorithm to produce a rational approximation for the leading-order solution based on the sampled data, $\hat{u}_0$ \eqref{e.LOapp}. We then performed an exponential asymptotic analysis and found that each pair of poles in the solution produced Stokes curves, each of which generated exponentially small oscillations. Taking the sum of these oscillations produced the total exponentially small behaviour $\hat{u}_{\mathrm{exp}}$ \eqref{e.AAAuexp}. 

Finally, we compared $u_{\mathrm{exp}}$ and $\hat{u}_{\mathrm{exp}}$, and found that $\hat{u}_{\mathrm{exp}}$ is able to accurately approximate $u_{\mathrm{exp}}$ for nonzero values of $\epsilon$, despite having different asymptotic behaviour in the limit that $\epsilon \to 0$. If $\epsilon$ is too small, however, the approximation is inaccurate. By reducing the tolerance of the AAA algorithm, and therefore increasing the accuracy of the rational approximation, we can reduce the threshold value of $\epsilon$ beyond which the approximation loses accuracy. 

Empirical tests suggest that this threshold value of $\epsilon$ is proportional to the distance between the branch point in $u_0$ and the nearest pole in $\hat{u}_0$. This is a measure of how accurately the rational approximation predicts the true singularity location. While the true branch point (or other singularity) in $u_0$ cannot typically be calculated in practice, it is possible to estimate the true location by studying the rate at which the poles in $\hat{u}_0$ accumulate \cite{Trefethen2021}.

There are two promising avenues available for extending this method to study nonlinear ordinary differential equations. The first is to determine the asymptotic corrections to the approximated exponentials that are caused by nonlinear interactions between nearby simple poles. The second is to use asymptotic arguments to determine the strength of singularities in the leading-order behaviour, and apply the AAA algorithm to a mapped version of the data that contains simple poles instead of other singularities. We briefly showed in Figure \ref{Fig:Nonlin} that this method can be used successfully for a nonlinear ordinary differential equation, and plan to explore this and related problems in future work.

This approach has connections to recent work \cite{Costin_2020,costin2022uniformization} where Pad\'e approximants are generated from numerical truncated power series data (as is typically available in Borel plane analyses \cite{crew2022resurgent}). These approximations generate accumulations of simple poles near branch points in the original Borel-transformed function. The authors use conformal mappings to convert these into poles, which may be more easily studied. It is likely that a variant of the the conformal map techniques developed in \cite{costin2022uniformization} could be applied to the AAA approximation to $u_0(z)$ developed in the present work, as well as more complicated leading-order behaviour.

Finally, it is important to determine how robust this method is to noisy input data. The effect of noise on numerical rational approximation has been studied in \cite{costin2022noise}, the authors study the effects of noise on conformal maps generated using Pad\'{e} approximation, and in \cite{wilber2022data}, in which the authors study the effects of noise on rational approximations generated using the AAA algorithm. It would be valuable to use similar methods to study the effect of noise on exponential asymptotic analyses.

\section{Acknowledgements}
CJL has been supported by Australian Research Council Discovery Project DP190101190. The authors would like to thank the Isaac Newton Institute (INI) for Mathematical Sciences for support and hospitality during the programme “Applicable resurgent asymptotics: towards a universal theory”, where part of the work on this paper was undertaken. The INI programme was supported by the EPSRC grant no EP/R014604/1. The authors would like to thank the Okinawa Institute of Science and Technology (OIST) for support and hospitality. Part of this research was conducted while visiting OIST through the Theoretical Sciences Visiting Program (TSVP). CJL would like to thank Nick L. Trefethen and John R. King for valuable discussions.

\bibliographystyle{plain}
\bibliography{reference2}

\end{document}